\documentclass[11pt]{amsart}
\usepackage{amsmath,amsthm, amscd, amssymb, amsfonts, mathrsfs,pb-diagram,pb-xy}
\usepackage[all]{xy}
\usepackage{enumitem}
\usepackage{mathrsfs}
\usepackage{srcltx}
\usepackage{hyperref}
\usepackage{multicol}
\usepackage[dvips, dvipsnames, usenames]{color}

\usepackage{ textcomp }

\allowdisplaybreaks

\newcommand{\superqa}[3]{{\bf A}_{#1}(#2|#3)}

\newcommand{\lu}{\mathcal{L}}
\newcommand{\luq}[1]{\lu_{#1}}

\newcommand{\ad}{\operatorname{ad}}

\newcommand{\co}{\operatorname{co}}

\newcommand{\gr}{\operatorname{gr}}

\newcommand{\he}{\operatorname{ht}}
\newcommand{\id}{\operatorname{id}}
\newcommand{\ord}{\operatorname{ord}}

\newcommand{\GK}{\operatorname{GKdim}}
\newcommand{\supp}{\operatorname{supp}}

\newcommand{\ku}{ \Bbbk}
\newcommand{\I}{\mathbb I}
\newcommand{\G}{\mathbb G}
\newcommand{\N}{\mathbb N}
\newcommand{\Z}{\mathbb Z}

\newcommand\zt{\Z^{\theta}}

\newcommand{\toba}{\mathscr{B}}
\newcommand{\wtoba}{\widetilde{\mathscr{B}}}
\newcommand{\htoba}{\widehat{\toba}}
\newcommand{\bq}{\mathfrak{q}}

\newcommand{\bp}{\mathfrak{p}}
\newcommand{\Hc}{{\mathcal H}}

\newcommand{\ydh}{{}^{H}_{H}\mathcal{YD}}
\newcommand{\ydkg}{{}^{\ku\Gamma}_{\ku\Gamma}\mathcal{YD}}
\newcommand{\ydkzt}{{}^{\ku\Z^{\theta}}_{\ku\Z^{\theta}}\mathcal{YD}}

\newcommand{\Ic}{\mathcal{I}}

\newcommand{\Pc}{\mathcal{P}}

\newcommand{\Zc}{\mathcal{Z}}
\newcommand{\hZ}{\widehat{\mathcal{Z}}}
\newcommand{\Jc}{\mathcal{J}}

\newcommand{\hOc}[1]{\underline{\widehat{\fO}^{#1}_+}}

\def\g{\mathfrak{g}}

\def\sli{\mathfrak{sl}}

\newcommand{\pre}[1]{\mathfrak{Pre}(#1)}
\newcommand{\prefd}[1]{\mathfrak{Pre}_{\operatorname{fGK}}(#1)}
\newcommand{\pregr}[1]{\mathfrak{Pre}^{\gr}(#1)}
\newcommand{\prefdgr}[1]{\mathfrak{Pre}_{\operatorname{fGK}}^{\gr}(#1)}

\newcommand{\qti}{\widetilde{q}}

\newcommand{\pf}{\begin{proof}}
\newcommand{\epf}{\end{proof}}

\newcommand{\GL}{\operatorname{GL}}
\newcommand{\fO}{\mathfrak O}
\newcommand{\wfO}{\underline{\mathfrak O}}
\newcommand{\wbeta}{\underline{\beta}}
\newcommand{\walpha}{\underline{\alpha}}
\newcommand{\wgamma}{\underline{\gamma}}

\newcommand{\n}{\mathfrak{n}}

\newcommand{\comp}[1]{\mathcal{P}_{c}(#1)}

\newcommand{\coefd}[2]{{d}_{(#1|#2)}}

\numberwithin{equation}{section}
\theoremstyle{plain}

\newtheorem{theorem}{Theorem}[section]
\newtheorem{lemma}[theorem]{Lemma}
\newtheorem{definition-theorem}[theorem]{Definition-Theorem}
\newtheorem{prop}[theorem]{Proposition}

\newtheorem{question}[theorem]{Question}

\theoremstyle{definition}
\newtheorem{definition}[theorem]{Definition}
\newtheorem{example}[theorem]{Example}

\newtheorem{remark}[theorem]{Remark}

\newcommand{\ot}{\otimes}

\newcommand{\ugamma}{\underline{\gamma}}

\newcommand{\ubeta}{\underline{\beta}}

\begin{document}
\noindent
\title[Posets of finite GKdim pre-Nichols algebras of diagonal type]
{Posets of finite GK-dimensional graded pre-Nichols algebras of diagonal type}

\author[Angiono]{Iv\'an Angiono}
\email{ivan.angiono@unc.edu.ar}
\author[Campagnolo]{Emiliano Campagnolo}
\email{emiliano.campagnolo@mi.unc.edu.ar}

\keywords{Hopf algebras, Nichols algebras, Gelfand-Kirillov dimension.
\\
MSC2020: 16T05, 16T20, 17B37, 17B62.}

\thanks{The work of both authors was partially supported by CONICET and Secyt (UNC)}

\begin{abstract}
We classify graded pre-Nichols algebras of diagonal type with finite Gelfand-Kirillov dimension. The characterization is made through an isomorphism of posets with the family of appropriate subsets of the set of positive roots coming from central extensions of Nichols algebras of diagonal type, generalizing the corresponding extensions for small quantum groups in de Concini-Kac-Procesi forms of quantum groups. 

On the way to achieving this result, we also classify graded quotients of algebras of functions of unipotent algebraic groups attached to semisimple Lie algebras.
\end{abstract}

\maketitle

\section{Introduction}

Let $\Bbbk$ be an algebraically closed field of characteristic zero. 
The so called quantized enveloping algebras $U_q(\g)$, where $q$ is a parameter and $\g$ is a semisimple Lie algebra, emerged after the works of Drinfeld and Jimbo as examples of non-commutative and non-cococommutative Hopf algebras over the field $\Bbbk(q)$: They were obtained by \emph{deforming} the structure of the corresponding enveloping algebra over $\Bbbk(q)$. When we consider the evaluation of the parameter $q$ in elements of $\Bbbk$, we get a Hopf algebra over $\Bbbk$, which behaves as $U(\g)$ in terms of representations if $q$ is not a root of unity.
In the nineties, de Concini, Kac and Procesi \cite{DKP,DP} studied the case in which $q$ is a root of unity of order $N$, under some mild conditions on $N$, and found a completely different story. To begin with, the centre of $U_q(\g)$ is larger than in the case where $q$ is not a root of unity, and it contains a Hopf subalgebra $\Zc_q$ which gives rise to an extension of Hopf algebras
\begin{align}\label{eq:intro-extension-DKP}
\Zc_q \hookrightarrow U_q(\g) \twoheadrightarrow \mathfrak{u}_q(\g).
\end{align}
Here $\mathfrak{u}_q(\g)$ is the Frobenius-Lusztig kernel, a finite-dimensional pointed Hopf algebra. The name comes from the evaluation of a different form of $U_q(\g)$ studied by Lusztig \cite{Lu-book,Lu-paper} in connection with the representation theory of algebraic groups in positive characteristic. He took an integral form generated by \emph{divided powers} of the generators and the algebra $\mathcal{U}_q(\g)$ obtained after evaluation in $q$ fits into an extension of Hopf algebras
\begin{align}\label{eq:intro-extension-Lusztig}
\mathfrak{u}_q(\g)\hookrightarrow \mathcal{U}_q(\g) \twoheadrightarrow U(\g).
\end{align}

Coming back to \eqref{eq:intro-extension-DKP}, $U_q(\g)$ is $\Z^{\theta}$-graded, where $\theta$ is the rank of $\g$, and has a triangular decomposition $U_q(\g)\simeq U_q^+(\g)\ot U_q^0(\g) \ot U_q^-(\g)$, where $U_q^0(\g)$ is the group algebra of a free abelian group in generators $K_i$, $1\le i\le \theta$, and $U_q^{\pm}(\g)$ has a PBW basis made of PBW generators $E_{\alpha}$, respectively $F_{\alpha}$, of degree $\pm \alpha\in\varDelta_+$, where $\varDelta$ is the set of roots of $\g$ viewed as a subset of $\Z^{\theta}$. In addition, $\Zc_q$ is the subalgebra generated by $E_{\alpha}^N$, $F_{\alpha}^N$, $\alpha\in\varDelta_{+}$, and $K_i^N$, $1\le i\le \theta$, while $\mathfrak{u}_q(\g)$ has a \emph{restricted} PBW basis with the same set of generators, but we restrict the powers up to $N$. Also, the restriction of \eqref{eq:intro-extension-DKP} to the corresponding positive parts gives an extension
$\Zc_q^+ \hookrightarrow U_q^+(\g) \twoheadrightarrow \mathfrak{u}_q^+(\g)$ of Hopf algebras, but in the braided tensor category $\ydkzt$ of Yetter-Drinfeld modules over $\zt$.

\medspace

A few years later, Andruskiewitsch and Schneider \cite{AS-lifting} introduced the so called Lifting Method to classify finite-dimensional pointed Hopf algebras.
In a nutshell, this method is based on the decomposition of the associated coradically graded Hopf algebra into the bosonization between a group algebra $\Bbbk\Gamma$ and a coradically graded Hopf algebra $R$ in the category $\ydkg$ by solving the following steps: first, to classify all finite dimensional Nichols algebras (see \S \ref{subsec:nichols-pre-N} and references therein for the definition and examples); then classify the possible finite dimensional Hopf algebras $R$ extending the Nichols algebras as before, called \emph{post-Nichols algebras}; and finally to obtain all the liftings of the corresponding bosonizations. This method was widely applied, with the first main result \cite{AS-Annals} being the classification in the case of abelian groups with moderate restriction on the order, where the examples are certain deformations of the bosonizations of $\mathfrak{u}_q^+(\g)$ with appropriate abelian groups. The general answer for abelian groups involves the classification of finite-dimensional Nichols algebras $\toba_{\bq}$ of diagonal type (depending on a braiding matrix $\bq\in(\Bbbk^{\times})^{\theta\times\theta}$), which in turn is contained in the classification of those with finite root system (i.e., a finite number of PBW generators) given by Heckenberger \cite{H-classif}. To do so, we attach a (kind of) Dynkin diagram with labels depending on $\bq$ and consider the connected components of this diagrams: the root system is, as expected, the disjoint union of the root systems of the connected components, so the list in \cite{H-classif} contains just those matrices $\bq$ with connected Dynkin diagram. In \cite{AA17} this list was split into the following families:
\begin{itemize}
\begin{multicols}{3}
\item Cartan type
\item Super type
\item Standard type
\item Modular type
\item Supermodular type
\item UFO.
\end{multicols}
\end{itemize}
Here, Cartan type is essentially the case of quantized enveloping algebras, while super type is related with quantized enveloping Lie superalgebras.

\medbreak

Andruskiewitsch and Schneider also started the classification of pointed Hopf algebras with finite Gelfand-Kirillov dimension (i.e. infinite-dimensional ones with some kind of moderate growth) in \cite{AS-crelle} by classifying those which are domains and satisfy a technical condition. To do so, they extend the Lifting Method to this context and obtained that all possible Nichols algebras are close to $U_q^+(\g)$ for $q$ not a root of unity, and the unique possible Hopf algebras $R$ extending the Nichols algeras are just the Nichols algebras themselves. 

If we want to get all Hopf algebras including those which are not domains, the answer is fully different: indeed, Lusztig examples provide post Nichols algebras $\mathcal{U}_q^+(\g)$ properly extending the Nichols algebras $\mathfrak{u}_q^+(\g)$. Taking graded duals, the extension between the positive parts in \eqref{eq:intro-extension-Lusztig} becomes that of \eqref{eq:intro-extension-DKP}, where $U_q^+(\g)$ is a pre-Nichols algebra (a graded intermediate quotient between the tensor algebra and the Nichols algebra $\mathfrak{u}_q^+(\g)$) of finite $\GK$, and classifying post-Nichols algebras with finite $\GK$ is related to classifying pre-Nichols algebras with finite $\GK$. In general, finite-dimensional Nichols algebras $\toba_{\bq}$ of diagonal type fit into an exact sequence of braided Hopf algebras $\Zc_{\bq}^+ \hookrightarrow \wtoba_{\bq} \twoheadrightarrow \toba_{\bq}$, generalizing the one between positive parts in \eqref{eq:intro-extension-DKP}, where $\wtoba_{\bq}$ is the distinguished pre-Nichols algebra \cite{An-dist}: it has a PBW basis with the same set of generators as $\toba_{\bq}$, but where we allow the powers of some of them to be arbitrary as in \eqref{eq:intro-extension-DKP}. Thus we may identify first those Nichols algebras of diagonal type with finite $\GK$ and then obtain all possible pre-Nichols algebras of finite $\GK$ covering each one of these Nichols algebras. For the first question the answer was given in \cite{angiono-garcia-2022}, following the conjecture made in \cite{AAH-memoirs}: a Nichols algebra of diagonal type has finite $\GK$ if and only if its root system is finite, i.e. it appears in the lists in \cite{H-classif}. Thus we have to classify all pre-Nichols algebras $\toba$ with finite $\GK$ covering the Nichols algebras $\toba_{\bq}$ with finite root system. 

Fixing a braiding matrix $\bq$, the corresponding set of pre-Nichols algebras form a poset whose maximal element is $\toba_{\bq}$ and we may wonder if there exists a minimal element between those of finite $\GK$, called the \emph{eminent} pre-Nichols algebra $\htoba_{\bq}$ in \cite{ASa}.  This is the case for all $\bq$ with connected diagrams up to Cartan types $A_{\theta}$, $D_{\theta}$ with label $q=-1$, thanks to \cite{ASa,ACSa,ACSa2,C}: in most cases eminent and distinguished pre-Nichols algebras coincide, i.e. $\htoba_{\bq}=\wtoba_{\bq}$. In any case, $\htoba_{\bq}$ fits into an exact sequence of braided Hopf algebras 
$\hZ_{\bq} \hookrightarrow \htoba_{\bq} \twoheadrightarrow \toba_{\bq}$, where $\hZ_{\bq}$ is an algebra of $q$-polynomials whose variables are homogeneous: we collect the $\N_0^{\theta}$-degrees of these variables in a set denoted by $\hOc{\bq}$, which is the set of positive roots of a classical root system when $\htoba_{\bq}=\wtoba_{\bq}$ by \cite{AAR-imrn}. As we will explain later, we are interested in subsets $B\subseteq \hOc{\bq}$ closed by sums: it means if $\alpha,\beta\in B$ are such that $\alpha+\beta\in\hOc{\bq}$, then $\alpha+\beta\in B$. Let $\comp{\bq}$ be the set of all subsets of $\hOc{\bq}$ closed by sums, which is a subposet of the poset of subsets of $\hOc{\bq}$.

Due to the results stated above, the determination of the poset of pre-Nichols algebras with finite $\GK$ is equivalent to the characterization of all intermediate quotients $\toba$ between $\htoba_{\bq}$ and $\toba_{\bq}$. The main result of this paper deals with this question:

\begin{theorem}\label{thm:main}
Let $\bq$ be a braiding matrix whose connected components are not of Cartan types $A_{\theta}$, $D_{\theta}$ with label $q=-1$ neither one dimensional with label $q=\pm 1$.
For each $\ubeta\in\hOc{\bq}$ let $z_{\beta}$ be a generator of $\hZ_{\bq}$ of degree $\ubeta$, and for each $B\in\comp{\bq}$ let
\begin{align*}
	\toba(\bq,B)& :=\htoba_{\bq} / \langle z_{\ubeta} | \beta\in \hOc{\bq}-B \rangle.
\end{align*}
Then each $\toba(\bq,B)$ is an $\N_0^{\theta}$-graded pre-Nichols algebra such that 
$\GK \toba(\bq, B)= |B|$. The assignment $B\mapsto \toba(\bq,B)$ gives an anti-isomorphism of posets between
\begin{itemize}
\item the set $\comp{\bq}$ of all subsets of $\hOc{\bq}$ closed by sums, and 
\item the set of all $\N_0^{\theta}$-graded pre-Nichols algebras of $\bq$ with finite $\GK$.
\end{itemize}
\end{theorem}

In addition, we prove that the poset in the non-connected case decomposes as the product of the corresponding posets of the connected components and  give a closed formula for the Hilbert series of each pre-Nichols algebra $\toba(\bq,B)$ as well as a PBW basis.

Going back through the steps of the Lifting Method, we classify by taking graded duals all post-Nichols algebras of diagonal type and finite $\GK$ (up to few exceptions on the connected components), which in turn give all coradically graded pointed Hopf algebras with abelian coradical and finite $\GK$ after bosonization with suitable abelian groups.

\medbreak

We can observe that any $\toba(\bq,B)$ fits into an exact sequence of the Nichols algebra $\bq$ by a $q$-central Hopf subalgebra, so we may ask:

\begin{question}
Are there examples of pre-Nichols algebras of finite $\GK$ which are not a ``central'' extension of braided Hopf algebras of the corresponding Nichols algebra?
\end{question}

Although the restriction to the $\N_0^{\theta}$-graded case may seem very strong, it has both a \emph{realisation-independence} reason and also a reduction to a problem with a closed answer: the general case depends strongly on the realisation of the braided vector space of diagonal type as Yetter-Drinfeld module, and a general answer may be somewhat unmanageable, see Remark \ref{rem:restriction-N0theta-graded}.

\medbreak

The organization of the paper is the following. First we recall several notions about Nichols algebras, distinguished and eminent pre-Nichols algebras when the braiding is of diagonal type; we also summarize known results about eminent pre-Nichols algebras and solve some questions on quotients of pre-Nichols algebras, extensions of Nichols algebras by central subalgebras and twist equivalence of braidings of diagonal type in \S \ref{sec:Nichols-and-related}. Motivated by the results in this section we consider quotients of the algebras of functions of unipotent algebraic groups which are the positive parts of semisimple ones; hence, in \S \ref{sec:unipotent-groups} we give the classification of these quotients in terms of subsets closed by sums of the set of positive roots of the associated semisimple Lie algebra. Finally we attack the determination of the poset of $\N_0^{\theta}$-graded pre-Nichols algebras with finite $\GK$ of a matrix $\bq$ such that $\toba_{\bq}$ has a finite root system (or equivalently, such that $\GK\toba_{\bq}<\infty$). Due to the results in \S \ref{sec:Nichols-and-related}, we can relate these quotients to those of the skew central Hopf subalgebra $\hZ_{\bq}$ (the subalgebra of coinvariants of the projection $\htoba_{\bq}\twoheadrightarrow \toba_{\bq}$), and also we can move to the case in which $\hZ_{\bq}$ is central. We attack first the connected case: we apply results in \S \ref{sec:unipotent-groups} to solve all the cases where $\htoba_{\bq}=\wtoba_{\bq}$, and compute explicitly the poset for the few exceptions where $\htoba_{\bq} \ne \wtoba_{\bq}$. Then we deal with the non-connected case using tools from \cite{ASa} and the answer for the connected case.

\subsection*{Notation}
We fix $\theta \in \N$ and set $\I = \I_{\theta} := \{1, 2, ..., \theta\}$.
Let $(\alpha_j)_{j\in \I}$ be the canonical basis of $\zt$ and $\alpha_{ij}$ denote $\alpha_{i}+\cdots+\alpha_j$, $i\leq j$.
Let $\beta =\sum_{i\in\I} a_i\alpha_i\in\Z^{\theta}$, sometimes also denoted $\beta=1^{a_1}\cdots \theta^{a_{\theta}}$ to shorten expressions.
The \emph{support} and the \emph{height} of $\beta$ are given by
\begin{align*}
\supp \beta &= \{i\in\I | a_i\neq 0 \}, & \he(\beta)=\sum_{i\in\I}a_i \in\Z.
\end{align*}
If $\gamma =\sum_{i\in\I} b_i\alpha_i\in\zt$ is such that $a_i\le b_i$ for all $i\in\I$, then we say that $\beta\le \gamma$.

If $N \in \N$ and $v\in \ku^{\times}$, then $(N)_v := \sum_{j=0}^{N-1}v^{j}$.
We denote by $\G_N$ to the group of roots of unity of order $N$ in $\Bbbk$, and $\G_N'$ the subset of primitive roots of order $N$.

Let $A$ be an associative algebra (with unit). We denote by $\GK A$ the Gelfand-Kirillov dimension of $A$. We refer to \cite{KL} for the definition and properties.

We will deal with $\N_0^{\theta}$-graded objects $U=\oplus_{\alpha\in \N_0^{\theta}}U_{\alpha}$. 
The Hilbert series of $U$ is
\begin{align*}
\Hc_{U} &= \sum_{\alpha\in\N_0^{\theta}} \dim U_{\alpha} \, t^{\alpha}\in \N_0[[t_1,\dots,t_{\theta}]]
 && \text{where }t^{\alpha}=t_1^{a_1}\cdots t_{\theta}^{a_{\theta}} \text{ for }\alpha=(a_1,\cdots, a_{\theta}).
\end{align*}
Given $\Hc=\sum_{\alpha\in\N_0^{\theta}} a_{\alpha} \, t^{\alpha},\Hc'=\sum_{\alpha\in\N_0^{\theta}} b_{\alpha} \, t^{\alpha}\in \N_0[[t_1,\dots,t_{\theta}]]$, we say that $\Hc\ge \Hc'$ if $a_{\alpha}\ge b_{\alpha}$ for all $\alpha\in\N_0^{\theta}$. Thus, if $U' \subseteq U$ as $\N_0^{\theta}$-graded objects, then $\Hc_U\ge \Hc_{U'}$.

Let $C$ be a coalgebra. We will use Sweedler notation for $C$ and any (left) comodule $V$; explicitly, $\Delta(c)=c_{(1)}\ot c_{(2)}\in C\ot C$ for all $c\in C$, and if $\rho:V\to C\ot V$ is the coaction, then $\rho(v)=v_{(-1)}\ot v_{(0)}$ for all $v\in V$.

\section{On pre-Nichols algebras of diagonal type}\label{sec:Nichols-and-related}

We start by recalling notions and results related with Nichols and pre-Nichols algebras, with special focus on the diagonal case.
Let $H$ be a Hopf algebra. As usual,  we denote by $\ydh$ the category of (left) Yetter-Drinfeld modules over $H$. We refer to \cite{A-leyva,Rad-libro} 
for unexplained notions and notations on Yetter-Drinfeld modules and braided vector spaces, and to \cite{AA17} for more information on Nichols algebras of diagonal type, and to \cite{KL} for definitions and basic results on Gelfand-Kirillov dimension.

\subsection{Nichols algebras and pre-Nichols algebras}\label{subsec:nichols-pre-N}

Recall that $\ydh$ is a braided tensor category: 
for each pair $V,W\in\ydh$, the braiding is given by
\begin{align*}
c_{V,W} &:V\ot W\to W\ot V, & c_{V,W}(x \ot y)&=x_{(-1)}\cdot y \ot x_{(0)}, & &x\in V,y\in W.
\end{align*}
Therefore each pair $(V,c_{V,V})$, $V\in\ydh$, is a braided vector space.

The tensor algebra $T(V)= \bigoplus_{n\geq 0} V^{\ot n}$ becomes a graded Hopf algebra in $\ydh$ by declaring that every element in $V$ is primitive. 
The Nichols algebra $\toba(V)$ of $V$ is the quotient of $T(V)$ by the maximal Hopf ideal $\Jc(V)=\oplus_{n\geq2}\Jc^n(V)$ generated by homogeneous elements of degree $\ge 2$. Hence, $\toba(V)$ is an $\N_0$-graded Hopf algebra over $\ydh$, where the degree one part is $V$, coincides with the set of primitive elements and generates $\toba(V)$.

It is known that the structure of the Nichols algebra $\toba(V)$ depends on the braiding $c:=c_{V,V}\in\GL(V^{\ot 2})$, not really on the realisation as a Yetter-Drinfeld module. This is why we consider braided vector spaces throughout this paper: i.e. pairs $(V,c)$, where $V$ is a vector space and $c\in\GL(V^{\ot 2})$ is a solution of the braid equation.

\smallbreak

Prominent examples are braided vector spaces of \emph{diagonal type}. It means that there exist a basis $\{x_i\}_{i\in \I}$ and a matrix $\bq=(q_{ij})\in(\Bbbk^{\times})^{\I\times\I}$ such that the braiding is
\begin{align*}
c^{\bq}: & V\ot V\to V\ot V, & c^{\bq}(x_i\ot x_j)&=q_{ij}\,x_j\otimes x_i, \, i,j\in\I.
\end{align*}
The matrix $\bq$ is called the braiding matrix. The information of $\bq$ is encoded in the associated Dynkin diagram. This is a labelled graph with $\theta$ vertices, each of them labelled with $q_{ii}$, and a edge between vertices $i$ and $j$ if $\qti_{ij}:=q_{ij}q_{ji} \ne 1$, labelled with this scalar. Different braiding matrices can have the same Dynkin diagram: the associated Nichols algebras are not isomorphic but equivalent in some sense as we will see in \S \ref{subsec:twist-equivalent}.

Nichols algebras of diagonal type depend only on $\bq$, so we will denote it by $\toba_{\bq}$; we denote accordingly $\Jc_{\bq}$ to the defining ideal of $\toba_{\bq}$. In addition, $(V,c)$ is realized as a Yetter-Drinfeld module over $\Bbbk\zt$ in a canonical way: We set the coaction on $V$ given by $\rho(x_i)=\alpha_i\ot x_i$ and the action given by $\alpha_i\cdot x_j=q_{ij} \, x_j$, $i,j\in\I$. From here we can deduce that $\toba_{\bq}$ is $\zt$-graded, where each $x_i$ has degree $\alpha_i$.

\smallbreak

A pre-Nichols algebra of $V$ is a braided Hopf algebra $\toba$ which is the quotient of braided Hopf algebras of $T(V)$  by an $\N_0$-homogeneous Hopf ideal $\Jc=\oplus_{n\geq2}\Jc^n$. Thus $\Jc\subseteq \Jc(V)$ and there exist canonical graded Hopf algebra epimorphisms
$$T(V)\twoheadrightarrow \toba \twoheadrightarrow \toba(V)$$
whose restriction to degree one is $\id_V$. 
The set of pre-Nichols algebras of $V$ becomes a poset $\pre{V}$, where $ \toba_1 \leq \toba_2$ if $\id_V$ induces an epimorphism of braided Hopf algebras $ \toba_1 \twoheadrightarrow \toba_2$. This poset has maximum and minimum elements, i.e. $\toba(V)$ and $T(V)$.

Assume that $\GK\toba(V)<\infty$. The subset $\prefd{V}$ of pre-Nichols algebras of $V$ with finite Gelfand-Kirillov dimension is a subposet with maximum element $\toba(V)$.
In case it admits a minimum $ \htoba(V)$, we will say that $ \htoba(V)$ is the \emph{eminent pre-Nichols algebra} of $V$.
The existence and computation of eminent pre-Nichols algebras $\htoba(V)$ reduces the problem of finding the set of all pre-Nichols algebras of $V$ with finite $\GK$ to the problem of finding quotients of $\htoba(V)$. As we will recall in \S \ref{subsec:dist-pre-Nichols}, eminent pre-Nichols algebras exist for most $V$ of diagonal type.

\subsection{Central extensions of braided Hopf algebras}\label{subsec:central-ext}

As we want to study posets of pre-Nichols algebras, we have to deal with extensions of connected Hopf algebras in $\ydh$. 
Motivated by \cite[Theorem 3.2]{Tk-quotient} and \cite[Proposition 3.6]{A+} we can state the following correspondence between Hopf ideals and 
normal coideal subalgebras in $\ydh$.

\begin{prop}\label{prop:HopfIdeal-Subcoalgebras}
Let $R$ be a connected Hopf algebra in $\ydh$. The assignements
\begin{align*}
&\mathtt{A} \mapsto \mathtt{I}(\mathtt{A}):=R/R\mathtt{A}^+, & & I\mapsto \mathtt{A}(I):={}^{\co R/I} R,
\end{align*}
give a bijection between the set of normal right coideal subalgebras $\mathtt{A}$ of $R$ and the set of Hopf ideals $I$ of $R$.
\end{prop}
\pf
By \cite[Proposition 3.6 (c)]{A+} we have a bijection between the set of right coideal subalgebras $\mathtt{A}$ of $R$ and the set of coideals $I$ of $R$ whose quotient map is of $R$-modules. Now $\mathtt{A}$ is normal if and only if $\mathtt{I}(\mathtt{A})$ is an ideal: the proof is by direct computation, analogous to that in \cite[Proposition 1.4]{Tk-quotient}.
\epf

Recall that  an \emph{extension of braided Hopf algebras} \cite[\S 2.5]{AN} is a sequence of morphisms of braided Hopf algebras
$\Bbbk \rightarrow \mathtt{A} \overset{\iota}{\to} \mathtt{C} \overset{\pi}{\to} \mathtt{B} \rightarrow \Bbbk $
such that $\iota$ is injective, $\pi$ is surjective, $\ker \pi = \mathtt{C}\iota(\mathtt{A}^+)$ and $\mathtt{A}=\mathtt{C}^{\,\co\pi}$.
For the sake of simplicity we just write $\mathtt{A} \overset{\iota}{\hookrightarrow} \mathtt{C} \overset{\pi}{\twoheadrightarrow} \mathtt{B}$.
If $\mathtt{C}$ is connected then any surjective braided Hopf algebra morphism $\mathtt{C} \overset{\pi}{\twoheadrightarrow} \mathtt{B}$ gives an extension by choosing $\mathtt{A}=\mathtt{C}^{\,\co\pi}$, see  \cite[3.6]{A+}.

We say that an extension is \emph{central} if $\mathtt{A}$ is contained in the center of $\mathtt{C}$.

In case that $\mathtt{A}$, $\mathtt{B}$, $\mathtt{C}$ are $\N_0^{\theta}$-graded with finite-dimensional homogeneous components and the maps $\iota$, $\pi$ preserve the $\N_0^{\theta}$-grading, the Hilbert series of these algebras satisfy the equality
$\Hc_{\mathtt{C}}= \Hc_{\mathtt{A}}\Hc_{\mathtt{B}}$, cf. \cite[Lemma 2.4]{ACSa}.

\smallbreak

Now we deal with central extensions of Hopf algebras in $\ydh$ whose right hand side term is a Nichols algebra.

\begin{prop}\label{prop:HopfidealB-HopfidealZ}
Let $\Zc {\hookrightarrow} \toba \overset{\pi}{\twoheadrightarrow} \toba(V)$ be a central extension of connected graded braided Hopf algebras, where $\toba=\oplus_{n\geq0}\toba^n$ is a pre-Nichols algebra of $V$, i.e. $\toba^1=V$.

\begin{enumerate}[leftmargin=*,label=\rm{(\roman*)}]
\item\label{item:HopfidealB-HopfidealZ-bijection} The assignment $I\mapsto \toba I$ is a bijective correspondence between graded Hopf ideals of $\Zc$ and graded Hopf ideals of $\toba$ generated in degree $\ge 2$.
\item\label{item:HopfidealB-HopfidealZ-graded}  Assume that $\toba$, $\toba(V)$ are $\N_0^{\theta}$-graded and $\pi$ preserves the $\N_0^{\theta}$-grading, so $\Zc$ is also $\N_0^{\theta}$-graded. 
Let $I$ be a $\N_0^{\theta}$-graded Hopf ideal of $\Zc$, $\toba':=\toba/\toba I$, $\mathtt{Z}=\Zc^{\,\co\Zc/I}$. Then $\Hc_{\toba}=\Hc_{\toba'}\Hc_{\mathtt{Z}}$.
\end{enumerate}\end{prop}
\pf
\ref{item:HopfidealB-HopfidealZ-bijection} 
Let $A=\oplus_{n\geq0} A^n$ be a graded right coideal subalgebra of $\toba$ such that $A^1=0$. We claim that $A^n\subset\Zc$ for all $n\ge 0$: The proof is by induction on $n$. For $n=0$, $A^0=\Bbbk1\subset \Zc$. Now assume that $A^k\subseteq \Zc$ for all $k\le n$. For each $x\in A^{n+1}$, 
\begin{align*}
\Delta(x) -x\ot 1-1\ot x \in \oplus_{i=1}^n A_i \ot \toba_{n+1-i}.
\end{align*}
Thus, by inductive hypothesis, $(\pi\ot\id)\Delta(x)=\pi(x)\ot 1+1\ot x$, so
\begin{align*}
\Delta(\pi(x))=(\pi\ot\pi)\Delta(x)=\pi(x)\ot 1+1\ot \pi(x).
\end{align*}
Hence $\pi(x)$ is a primitive element of $\toba(V)$ in degree $n+1\ge 2$, so $\pi(x)=0$. Thus $x\in{}^{\co \pi}\toba=\Zc$, and the inductive step follows.
Therefore, the set of graded right coideal subalgebras of $\toba$ such that $A^1=0$ is the set of graded coideal subalgebras of $\Zc$, and all these coideal subalgebras are normal since $\Zc$ is central.

Using Proposition \ref{prop:HopfIdeal-Subcoalgebras}, we obtain a bijective correspondence between the corresponding Hopf ideals; that is, between the graded Hopf ideals of $\Zc$ and the graded Hopf ideals of $\toba$ generated in degree $\ge 2$.

\noindent\ref{item:HopfidealB-HopfidealZ-graded} Notice that $\Hc_{\Zc}=\Hc_{\Zc'}\Hc_{\mathtt{Z}}$ and $\Zc' {\hookrightarrow} \toba' \overset{\pi}{\twoheadrightarrow} \toba(V)$ is also a central extension of connected graded braided Hopf algebras. From this exact sequence and the one of $\toba$ we have that 
$\Hc_{\toba}= \Hc_{\Zc}\Hc_{\toba(V)}$ and $\Hc_{\toba'}= \Hc_{\Zc/I}\Hc_{\toba(V)}$. Thus the statement follows by putting together the three equalities involving Hilbert series. 
\epf

\subsection{Distinguished pre-Nichols algebras}\label{subsec:dist-pre-Nichols}

Let $(V,c^{\bq})$ be a braided vector space of diagonal type with $\GK\toba_{\bq}<\infty$.
By \cite{angiono-garcia-2022}, this means that $\toba_{\bq}$ has a PBW basis with a finite set of homogeneous generators; in other words, $\bq$ belongs to the lists in \cite{H-classif}.
The set $\varDelta_+^{\bq}$ of \emph{positive roots} of $\bq$ consist of the $\N_0^{\theta}$-degrees of these generators, which is independent of the chosen PBW basis \cite{AA17,H-Weyl gpd}. The \emph{set of roots} of $\bq$ is $\varDelta^{\bq}:=\varDelta_{+}^{\bq}\cup (-\varDelta_{+}^{\bq})$.

Assume from now on that $|\varDelta_+^{\bq}|<\infty$. For $\alpha=(a_1,\cdots,a_{\theta}),\beta=(b_1,\cdots,b_{\theta})\in\N_0^{\theta}$ set
\begin{align*}
q_{\alpha\beta}&:=\prod_{i,j=1}^{\theta}q_{ij}^{a_ib_j}, & N_{\beta}:=\ord q_{\beta\beta}\in\N\cup\{\infty\}.
\end{align*}
A total order $\succ$ on $\varDelta_+^{\bq}$ is \emph{convex} if for all $\alpha \succ \beta  \in \varDelta_+^{\bq}$ such that $\alpha+\beta\in \varDelta_+^{\bq}$ we have that $\alpha\succ\alpha+\beta\succ\beta$. For each convex order $\beta_1\succ \cdots \succ \beta_M$ there exists a PBW basis with set of generators $x_{\beta}$, $\beta\in\varDelta_{+}^{\bq}$, with each $x_{\beta}$ of degree $\beta$. More explicitly, the set
\begin{align*}
&x_{\beta_1}^{n_1} \cdots x_{\beta_M}^{n_M}, && 0\le n_i <N_{\beta_i},
\end{align*}
is a basis of $\toba_{\bq}$, see e.g. \cite{An-jems,HY-groupoid}. Thus, the Hilbert series of $\toba_{\bq}$ is
\begin{align*}
\Hc_{\toba_{\bq}}(t) =\left(\prod_{\beta\in\varDelta_{+}^{\bq}:N_{\beta}=\infty}\frac{1}{1-t^{\beta}} \right)
\left(\prod_{\beta\in\varDelta_{+}^{\bq}:N_{\beta}<\infty}\frac{1-t^{N_{\beta}\beta}}{1-t^{\beta}} \right).
\end{align*}

Next we move to pre-Nichols algebras of diagonal type and connected Dynkin diagram.
Between the Nichols algebras $\toba_{\bq}$ with finite $\GK$, some of them are infinite-dimensional. By \cite{C} for most of these $\bq$ the Nichols algebra $\toba_{\bq}$ is the unique pre-Nichols algebra with finite $\GK$, and for the remaining ones there exist exactly one proper pre-Nichols algebra $\htoba_{\bq}$ with finite $\GK$ (which is then eminent).

Thus we can restrict to the problem of determining pre-Nichols algebras of finite $\GK$ when $\dim\toba_{\bq}<\infty$. This is equivalent to the fact that $N_{\beta}<\infty$ for all $\beta\in\varDelta_{+}^{\bq}$. In this case, there exists a pre-Nichols algebra with finite $\GK$, called the \emph{distinguished pre-Nichols algebra} $\wtoba_{\bq}$ \cite{An-dist}. This pre-Nichols algebra is the quotient $\wtoba_{\bq}=T(V)/\Ic_{\bq}$ by the ideal $\Ic_{\bq}$ generated by the set of defining relations in \cite[Theorem 3.1]{An-crelle} adding a few extra relations and removing relations of the form $x_{\alpha}^{N_{\alpha}}$ for $\alpha\in\fO^{\bq}_+$. Here,
\begin{align*}
\fO^{\bq} :=\left\{ \alpha\in\varDelta^{\bq}: q_{\alpha\beta}q_{\beta\alpha}\in\{q_{\alpha\alpha}^n:n\in\Z\} \text{ for all }\beta\in\zt \right\}
\end{align*}
is the set of roots of Cartan type, cf. \cite{An-crelle,AAR-imrn}, and $\fO^{\bq}_+:=\fO^{\bq}\cap\N_0^{\theta}$. By \cite{An-dist} the set
\begin{align*}
&x_{\beta_1}^{n_1} \cdots x_{\beta_M}^{n_M}, && 0\le n_i <\widetilde{N}_{\beta_i},
\end{align*}
is a basis of $\wtoba_{\bq}$, where $\widetilde{N}_{\beta}:=\begin{cases} \infty, & \beta\in \fO^{\bq}_+, \\ N_{\beta}, & \beta\notin \fO^{\bq}_+.
\end{cases}$ Thus, the Hilbert series of $\wtoba_{\bq}$ is
\begin{align*}
\Hc_{\wtoba_{\bq}}(t) =\left(\prod_{\beta\in\fO^{\bq}_+}\frac{1}{1-t^{\beta}} \right)
\left(\prod_{\beta\in\varDelta_{+}^{\bq}-\fO^{\bq}_+}\frac{1-t^{N_{\beta}\beta}}{1-t^{\beta}} \right).
\end{align*}

Let $\Zc_{\bq}$ be the subalgebra of $\wtoba_{\bq}$ generated by $z_{\alpha}:=x_{\alpha}^{N_{\alpha}}$, $\alpha\in\fO^{\bq}_+$. By \cite{An-dist}
there exists an extension of braided Hopf algebras
\begin{align}\label{eq:central-ext-distinguished}
\Zc_{\bq} {\hookrightarrow} \wtoba_{\bq} \overset{\pi}{\twoheadrightarrow} \toba_{\bq},
\end{align}
i.e. $\Zc_{\bq}=\wtoba_{\bq}^{\co \pi}$. Moreover $\Zc_{\bq}$ is a $q$-central Hopf subalgebra of $\wtoba_{\bq}$, which is a $q$-polynomial algebra in variables $z_{\alpha}$, $\alpha\in\fO^{\bq}_+$. Now set
\begin{align*}
\wbeta&:=N_{\beta} \beta, \, \beta\in\fO^{\bq}, & \wfO^{\bq} &:=\{\wbeta: \beta\in\fO^{\bq}\}, & \wfO^{\bq}_+&:=\fO^{\bq}\cap\N_0^{\theta}.
\end{align*}
By \cite[Theorem 3.7]{AAR-imrn}, $\wfO^{\bq}$ is a root system (in the classical sense), with basis
\begin{align*}
\varPi^{\bq}:=\{ \wgamma\in \wfO^{\bq}_+: \wgamma\ne\walpha+\wbeta \text{ for all }\walpha,\wbeta\in\wfO^{\bq}_+ \}.
\end{align*}
With the notation above, the Hilbert series of $\wtoba_{\bq}$ can be also written as
\begin{align}\label{eq:Hilbert-series-distinguished}
\Hc_{\wtoba_{\bq}}(t) &= \Hc_{\toba_{\bq}}(t)\Hc_{\Zc_{\bq}}(t) = \Hc_{\toba_{\bq}}(t)\left(\prod_{\wbeta\in\wfO^{\bq}_+}\frac{1}{1-t^{\wbeta}} \right).
\end{align}

\begin{example}\label{ex:positive-part-Uq}
Fix $A=(a_{ij})$ a finite Cartan matrix, $(d_i)\in\N^{\theta}$ minimal such that $(d_ia_{ij})$ is symmetric and $q\in\Bbbk$ is a root of unity of order $N$ coprime with all $a_{ij}$'s.
Let $\bq=(q_{ij})$, where $q_{ij}=q^{d_ia_{ij}}$. In this case $\bq$ is of \emph{Cartan type} and $\wtoba_{\bq}\simeq U_q^+(\g)$, where $\g$ is the (finite-dimensional) semisimple Lie algebra with $A$ as Cartan matrix. Moreover
\begin{align*}
N_{\beta}&=N \text{ for all }\beta\in \Delta^{\bq}_+,  & \fO^{\bq}_+&=\Delta^{\bq}_+,  & \text{thus }\underline{\fO}^{\bq}_+&=\{N\beta: \beta\in \Delta^{\bq}_+\}.
\end{align*}

By \cite{DKP} $\Zc_{\bq}$ is the algebra of functions of the unipotent algebraic group with Lie algebra $\n_+$ (the positive part of $\g$, for a fixed Borel subalgebra).
\end{example}

\medbreak

From the Hilbert series we check that $\GK\wtoba_{\bq}=|\fO_+^{\bq}|<\infty$. Thus we may wonder if $\wtoba_{\bq}$ is the eminent pre-Nichols algebra of $\bq$. This is mostly the case. More precisely:

\begin{theorem}\label{thm:eminent-distinguished}\cite{ASa,ACSa,ACSa2}
Let $(V,\bq)$ be a braided vector space of diagonal type such that $\dim \toba_{\bq}<\infty$ and the Dynkin diagram is connected. Then the distinguished pre-Nichols algebra $\wtoba_{\bq}$ is eminent, except in the following cases:
\begin{enumerate}[leftmargin=*,label=\rm{(\Alph*)}]
\item\label{item:eminent-exception-AD} Cartan $A_{\theta}$ or $D_{\theta}$ with $q=-1$,
\item\label{item:eminent-exception-A2} $A_2$ with $q\in \G_3'$,
\item\label{item:eminent-exception-superA} $\superqa{3}{q}{\{2\}}$ or $\superqa{3}{q}{\{1,2,3\}}$, with $q\in \G_{\infty}$,
\item\label{item:eminent-exception-g(2,3)} $\g(2,3)$ with any of the following Dynkin diagram
\begin{align*}
\xymatrix@C=30pt{\overset{-1}{\circ} \ar@{-}[r]^{\xi} &\overset{-1}{\circ}\ar@{-}[r]^{\xi}\ &\overset{-1}{\circ}},&
&\xymatrix@C=30pt{\overset{-1}{\circ} \ar@{-}[r]^{\xi^2} &\overset{\xi}{\circ}\ar@{-}[r]^{\xi}\ &\overset{-1}{\circ}}.
\end{align*}
\end{enumerate}
\end{theorem}

If $\bq$ is as in \ref{item:eminent-exception-AD}, it is not even known whether the eminent pre-Nichols algebra exists. But for the other cases, there is an answer: the eminent pre-Nichols algebra is a $q$-central extension of the distinguished pre-Nichols algebra, as we will describe below.

\begin{theorem}\label{thm:eminent-not-distinguished}
Let $(V,\bq)$ be a braided vector space of diagonal type such that $\dim \toba_{\bq}<\infty$. 
\begin{enumerate}[leftmargin=*,label=\rm{(\alph*)}]
\item\label{item:eminent-not-distinguished-A2} \cite{ASa} If $\bq$ is of type $A_2$ with $q\in \G_3'$, then the eminent pre-Nichols algebra of $\bq$ is
\begin{align*}
\htoba_{\bq}=T(V)/\langle x_{1112}, x_{2112}, x_{2221}, x_{1221} \rangle,
\end{align*}
and $\GK \htoba_{\bq}=5$.
Let $Z_{\bq}$ be the subalgebra of $\htoba_\bq$ generated by $x_{112}$, $x_{221}$.
There is a $\N_0^2$-homogeneous $q$-central extension of braided Hopf algebras $Z_{\bq} \hookrightarrow \htoba_{\bq} \twoheadrightarrow \wtoba_{\bq}$, and $Z_{\bq}$ is a $q$-polynomial algebra in variables $x_{112}$ and $x_{221}$.

\item\label{item:eminent-not-distinguished-Asuper-1} \cite{ACSa} If $\bq$ is of type $\superqa{3}{q}{\{2\}}$ with $q\in \G_{N}'$, then 
\begin{align*}
\htoba_{\bq}=T(V)/\langle x_2^2,x_{13},x_{112},x_{332} \rangle
\end{align*}
is the eminent pre-Nichols algebra of $\bq$, and $\GK \htoba_{\bq}= 3$.
Let $Z_{\bq}$ be the subalgebra of $\htoba_\bq$ generated by $x_u:=[x_{123},x_2]_c$.
There is a $\N_0^3$-homogeneous $q$-central extension of braided Hopf algebras $Z_{\bq} \hookrightarrow \htoba_{\bq} \twoheadrightarrow \wtoba_{\bq}$, and $Z_{\bq}$ is a $q$-polynomial algebra in $x_u$.

\item\label{item:eminent-not-distinguished-Asuper-2} \cite{ACSa} If $\bq$ is of type $\superqa{3}{q}{\{1,2,3\}}$, with $q\in \G_{N}'$, then 
\begin{align*}
\htoba_{\bq}=T(V)/\langle x_1^2, x_2^2, x_3^2, x_{213}, [x_{123},x_2]_c \rangle
\end{align*}
is the eminent pre-Nichols algebra of $\bq$, and $\GK \htoba_{\bq}= 3$.
Let $Z_{\bq}$ be the subalgebra of $\htoba_\bq$ generated by $x_{13}$.
There is a $\N_0^3$-homogeneous $q$-central extension of braided Hopf algebras $Z_{\bq} \hookrightarrow \htoba_{\bq} \twoheadrightarrow \wtoba_{\bq}$, and $Z_{\bq}$ is a $q$-polynomial algebra in $x_{13}$.

\item\label{item:eminent-not-distinguished-g(2,3)-1} \cite{ACSa2} If $\bq$ is of type $\g(2,3)$ with diagram $\xymatrix@C=30pt{\overset{-1}{\circ} \ar@{-}[r]^{\xi} &\overset{-1}{\circ}\ar@{-}[r]^{\xi}\ &\overset{-1}{\circ}}$, then set
\begin{align*}
	x_u&:=[[x_{12},x_{123}]_c,x_2]_c, & x_v&:=[[x_{123},x_{23}]_c,x_2]_c,
\end{align*}
Then the eminent pre-Nichols algebra is
\begin{align*}
\htoba_{\bq}=&T(V)/\langle x_{1}^2,x_2^2,x_3^2,x_{13}, [x_1,x_u]_c,[x_1,x_v]_c,[x_u,x_3]_c,[x_v,x_3]_c\rangle,
\end{align*}
and $\GK \htoba_{\bq} =6$. Let $Z_{\bq}$ be the subalgebra of $\htoba_\bq$ generated by $x_u$ and $x_v$.
There is a $\N_0^3$-homogeneous $q$-central extension of braided Hopf algebras $Z_{\bq} \hookrightarrow \htoba_{\bq} \twoheadrightarrow \wtoba_{\bq}$, and $Z_{\bq}$ is a $q$-polynomial algebra in variables $x_u$, $x_v$.

\item\label{item:eminent-not-distinguished-g(2,3)-2} \cite{ACSa2} If $\bq$ is of type $\g(2,3)$ with diagram
$\xymatrix@C=30pt{\overset{-1}{\circ} \ar@{-}[r]^{\xi^2} &\overset{\xi}{\circ}\ar@{-}[r]^{\xi}\ &\overset{-1}{\circ}}$, then set 
\begin{align*}
	x_u&:=[[x_{123},x_2]_c,x_2]_c, & x_v&:=[x_{123},x_{12^23}]_c,
\end{align*}
Then the eminent pre-Nichols algebra is
\begin{align*}
	\htoba_{\bq}=&T(V)/\langle x_{1}^2,x_3^2,x_{13}, [x_{223},x_{23}]_c, x_{221},x_{2223}, [x_v,x_3]_c, [x_{12^33^2},x_2]_c, [x_{12^33^2},x_3]_c\rangle.
\end{align*}
and $\GK \htoba_{\bq} =6$. Let $Z_{\bq}$ be the subalgebra of $\htoba_\bq$ generated by $x_u$ and $x_v$.
There is a $\N_0^3$-homogeneous $q$-central extension of braided Hopf algebras $Z_{\bq} \hookrightarrow \htoba_{\bq} \twoheadrightarrow \wtoba_{\bq}$, and $Z_{\bq}$ is a $q$-polynomial algebra in variables $x_u$, $x_v$.
\end{enumerate}
\end{theorem}

\medspace

We denote by $\pre{\bq}$ and $\prefd{\bq}$ the corresponding posets of pre-Nichols algebras when $V$ is of diagonal type with matrix $\bq$.
Let $\pregr{\bq}$ be the subposet of those pre-Nichols algebras which are $\N_0^{\theta}$-graded, and $\prefdgr{\bq}=\prefd{\bq}\cap\pregr{\bq}$, i.e. the subposet of those pre-Nichols algebras with finite $\GK$ which are $\N_0^{\theta}$-graded. The main result of this work is the characterization of $\prefdgr{\bq}$ for all those cases where the eminent pre-Nichols algebra of all connected components of $\bq$ is known.

\begin{remark}\label{rem:restriction-N0theta-graded}
The main reason behind the restriction to the subposet of $\N_0^{\theta}$-graded pre-Nichols algebras is that this is the set of all pre-Nichols algebras that can be realised in the category of Yetter-Drinfeld modules for any principal realization of $V$ over a group $\Gamma$. Recall that a principal realization of a braided vector space of diagonal type means that there exists a basis $(x_i)$ of $V$, elements $g_i\in Z(\Gamma)$ and $\chi_i\in\widehat{\Gamma}$ such that the coaction of $x_i$ is given by $g_i$ and $g\cdot x_i=\chi_i(g)x_i$ for all $g\in \Gamma$, so $q_{ij}=\chi_j(g_i)$ for all $i,j\in\I$.

For example, let $\bq=(q_{ij})$ be such that $\qti_{ij}=1$, $q_{ii}\in\G_{N_i}'$, for all $i\ne j\in\I$, where $N_i\in \N_{0}$. The distinguished pre-Nichols algebra is the so called quantum plane, $\htoba_{\bq} = T(V)/ \langle x_{ij} | i<j\in\I \rangle$. 
Fix $i\ne j$ such that $N_i=N_j$ and $a\in \Bbbk^{\times}$. As $x_i^{N_i},x_j^{N_j}$ are primitive elements of the same degree, the quotient
\begin{align*}
\toba =T(V)/\langle x_{k\ell }, \, k<\ell \in\I; \, x_i^{N_i}-a\, x_j^{N_j} \rangle= \htoba_{\bq}/\langle x_i^{N_i}-a\, x_j^{N_j} \rangle
\end{align*}
is a pre-Nichols algebra of $\bq$ of finite $\GK$. 
Fix also a principal realization over a group $\Gamma$: If either $g_i^{N_i}\ne g_j^{N_j}$ or else $\chi_i^{N_i}\ne \chi_j^{N_j}$, then $\toba$ is not an object in $\ydkg$.
\end{remark}

\subsection{Twist equivalence and pre-Nichols algebras}\label{subsec:twist-equivalent}

Recall that two braiding matrices $\bq=(q_{ij})_{i,j\in\I}$ and $\bp=(p_{ij})_{i,j\in\I}$ are \emph{twist-equivalent} if $q_{ii}=p_{ii}$ and $q_{ij}q_{ji}=p_{ij}p_{ji}$ for all $i\ne j\in\I$, cf. \cite[Definition 3.8]{AS-MSRI}. We write $\bq \sim\bp$. Notice that two matrices are twist-equivalent if and only if their Dynkin diagrams coincide.

Let $(V,c)$ and $(W,d)$ be the braided vector spaces of diagonal type with matrices $\bq$, $\bp$, respectively.
We want to relate Hopf ideals of $T(V)$ and $T(W)$. 

\begin{prop}\label{prop:twist-equiv-poset-preNichols}
Let $\bq=(q_{ij})_{i,j\in\I}$ and $\bp=(p_{ij})_{i,j\in\I}$ be twist-equivalent matrices.
\begin{enumerate}[leftmargin=*,label=\rm{(\roman*)}]
\item\label{item:twist-equiv-poset-preNichols-1} There exists an isomorphism of posets $\Psi:\pregr{\bq}\to\pregr{\bp}$ preserving the Hilbert series.
\item\label{item:twist-equiv-poset-preNichols-2} $\Psi$ restricts to an isomorphism $\Psi:\prefdgr{\bq}\to\prefdgr{\bp}$.
\end{enumerate}
\end{prop}

\pf

\ref{item:twist-equiv-poset-preNichols-1}
In \cite[Proposition 3.9]{AS-MSRI} the authors introduce a linear isomorphism $\psi:\toba_{\bq}\to\toba_{\bp}$, which is a coalgebra isomorphism: Let us recall more details about this isomorphism.
In loc. cit. the authors take the group cocycle
\begin{align*}
\sigma &: \zt\times\zt\to\Bbbk^{\times}, & \sigma(g_i,g_j) &=\begin{cases} p_{ij}q_{ij}^{-1}, & i\le j \\ 1 & i>j. \end{cases}
\end{align*}
For any Hopf algebra $R\in\ydkzt$, this group cocycle induces, up to projection, a Hopf cocycle $\sigma:R\#\Bbbk\zt \otimes R\#\Bbbk\zt\to \Bbbk$ and we  consider the Hopf algebra $(R\#\Bbbk\zt)_{\sigma}$: The coalgebra structure does not change and the canonical inclusion 
$\Bbbk\zt\hookrightarrow (R\#\Bbbk\zt)_{\sigma}$ and projection $(R\#\Bbbk\zt)_{\sigma}\twoheadrightarrow\Bbbk\zt$ are still Hopf algebra maps. Thus $(R\#\Bbbk\zt)_{\sigma}$ decomposes as $(R\#\Bbbk\zt)_{\sigma}\simeq R_{\sigma}\#\Bbbk\zt$. 
As stated in \cite[Lemma 2.12]{AS-MSRI}, the assignment $R\mapsto R_{\sigma}$ takes Hopf algebras in $\ydkzt$ to Hopf algebras in $\ydkzt$, where the coalgebra structure keeps unchanged, so it restricts to graded Hopf algebras: If $R=\oplus_{n\geq0} R^n$, then $R_{\sigma}=\oplus_{n\geq0} R_{\sigma}^n$, with $R_{\sigma}^n=R^n$ as vector spaces.

Let $V$, $W$ be as above: as usual we consider $V,W\in\ydkzt$. Working as in \cite[Proposition 3.9 \& Remark 3.10]{AS-MSRI} we check that $T(V)_{\sigma}=T(W)$. Any $\zt$-graded pre-Nichols algebra $\toba$ of $V$ is a Hopf alegebra in $\ydkzt$: the Hopf algebra projection $\pi:T(V)\twoheadrightarrow \toba$ gives a Hopf algebra projection $\pi:T(W)=T(V)_{\sigma}\twoheadrightarrow \toba_{\sigma}$, which preserves the $\N_0$-graded components. Thus, we have a map $\Psi:\pregr{\bq}\to\pregr{\bp}$, $\Psi(\toba)=\toba_{\sigma}$. Moreover, $\Psi$ is a map of posets, which has an inverse map given by $\sigma^{-1}$.

\medbreak

\noindent \ref{item:twist-equiv-poset-preNichols-2} As $\Psi$ preserves the Hilbert series, $\GK \toba_{\sigma}=\GK \toba$ by \cite[Lemma 6.1]{KL}.
\epf

Next we want to reduce to the case in which $\Zc_{\bp}$ is central (more than skew central). We can do this reduction up to twist-equivalence.

\begin{lemma}\label{lem:twist-equivalent-central}
Let $\bq$ be a matrix such that $\dim\toba_{\bq}<\infty$. There exists $\bp\sim\bq$ such that $\Zc_{\bp}$ is a central subalgebra of $\wtoba_{\bp}$.
\end{lemma}
\pf
By \cite[Remark 4.4]{AAY} we need to find a matrix $\bp\sim\bq$ such that $p_{\alpha_i\beta}^{N_{\beta}}=1$ for all $i\in\I$ and all $\beta\in\varPi_{\bp}$, and it is enough to check it just for one matrix $\bq$ in each Weyl equivalence class.

If $\bq$ belongs to the one-parameter families, then the existence of $\bp$ follows by \cite[Appendix A]{AAY}. The remaining cases are treated case-by-case.
\epf

Putting together Proposition \ref{prop:twist-equiv-poset-preNichols} and Lemma \ref{lem:twist-equivalent-central} we can restrict to matrices $\bq$ such that $\Zc_{\bq}$ is a central Hopf subalgebra. If so, then $\Zc_{\bq}$ is the algebra of functions of an algebraic group. Taking graded duals in the central extension in \eqref{eq:central-ext-distinguished} we get a new extension of braided graded connected Hopf algebras
\begin{align}\label{eq:central-ext-lusztig}
\toba_{\bq^{t}} {\hookrightarrow} \luq{\bq} \overset{\pi}{\twoheadrightarrow} \Zc_{\bq}^d.
\end{align}
Here $\luq{\bq}=\wtoba_{\bq^t}^d$ is called the \emph{Lusztig algebra} of $\bq$, cf. \cite{AAR-imrn}. Notice that $\Zc_{\bq}^d$ is cocommutative, thus it is (isomorphic to) the enveloping algebra of a finite-dimensional nilpotent Lie algebra $\n_{\bq}$: By \cite[Theorem 1.1]{AAR-imrn}, $\n_{\bq}$ is the positive part of the semisimple Lie algebra with root system $\wfO^{\bq}$.

\section{Hopf ideals of algebras of functions of unipotent groups}\label{sec:unipotent-groups}

According with the previous section we have to understand the coalgebra structure of the algebra of functions of the unipotent algebraic group whose Lie algebra is the positive part $\n_+$ of a semisimple Lie algebra $\g$ with root system $\varDelta$, that we denote $\Zc_{\varDelta}$. Let $\vartheta$ be the rank of $\varDelta$, $\I:=\{1,\cdots,\vartheta\}$.

As an algebra, $\Zc_{\varDelta}=\Bbbk[z_{\alpha}|\alpha\in\varDelta_+]$. Thus the set $M(\Zc_{\varDelta})$ of monomials in variables $z_{\alpha}$, $\alpha\in\varDelta_+$, is a basis of $\Zc_{\varDelta}$; we will say that a monomial $m\in M(\Zc_{\varDelta})$ \emph{appears} in an element $w\in\Zc_{\varDelta}$ if the coefficient of $m$ in the expression of $w$ as a linear combination of the elements of $M(\Zc_{\varDelta})$ is not zero. We also fix $M(\Zc_{\varDelta})\ot M(\Zc_{\varDelta})$ and $M(\Zc_{\varDelta})\ot M(\Zc_{\varDelta})\ot M(\Zc_{\varDelta})$ as bases of $\Zc_{\varDelta}\ot \Zc_{\varDelta}$ and $\Zc_{\varDelta}\ot \Zc_{\varDelta}\ot \Zc_{\varDelta}$, and say that $m\ot m'$, respectively $m\ot m'\ot m''$, appears in a term, correspondingly.

Let $Q_+$ be the lattice of positive roots, i.e. if $\{\alpha_i| i\in\I\}$ is the set of simple roots of $\varDelta_+$, then 
$Q_+=\sum_{i=1}^{\vartheta} \Z\alpha_i$. Then $\Zc_{\varDelta}$ is a $Q_+$-graded Hopf algebra, with each $z_{\alpha}$ of degree $\alpha$.

Now we relate $Z_{\varDelta}$ and $U(\n_+)$. Although this relation may be well known, we give a proof since the involved tools are useful in the forthcoming proofs.

\begin{lemma}\label{lem:structure-of-graded-dual-of-Z-Delta}
\begin{enumerate}[leftmargin=*]
\item\label{item:Z-Delta-graded-dual} The $Q_+$-graded dual of $\Zc_{\varDelta}$ is (isomorphic to) $U(\n_+)$, the enveloping algebra of the nilpotent Lie algebra $\n_+$.
\item\label{item:Z-Delta-primitive} The subspace $\Pc(\Zc_{\varDelta})$ of primitive elements of $\Zc_{\varDelta}$ has basis $\{z_{\alpha_i} \, |\, i\in\I\}$.
\end{enumerate}
\end{lemma}
\pf
For \eqref{item:Z-Delta-graded-dual}, let $U=\oplus_{\beta\in Q_+} U_{\beta}$ be the graded dual of $\Zc_{\varDelta}$, with $U_{\beta}$ the component of degree $\beta$. As $\Zc_{\varDelta}$ is connected (as algebra) and commutative, $U$ is connected (as coalgebra), and cocommutative. Thus $U$ is the enveloping algebra of its primitive elements.

The subspace $\Pc(U)$ of primitive elements of $U$ is the space of derivations $\Zc_{\varDelta}\to\Bbbk$, which is canonically identified with $\n_+$.

\smallskip
For \eqref{item:Z-Delta-primitive}, we proceed dually to the previous statement: $\Pc(\Zc_{\varDelta})$ is the space of derivations $\partial:U\to\Bbbk$.
This space has dimension $\le\vartheta$, since $U$ is generated by $e_i$, $i\in\I$, as algebra, so any derivation $\partial:U\to\Bbbk$ is univocally determined by the values $\partial(e_i)$, $i\in\I$. On the other hand, each $z_{\alpha_i}$ is primitive, so $\dim\Pc(\Zc_{\varDelta})\ge \vartheta$. From these two statements, $\dim\Pc(\Zc_{\varDelta})=\vartheta$, and $\Pc(\Zc_{\varDelta})$ has basis $\{z_{\alpha_i} \, |\, i\in\I\}$.
\epf

\begin{remark}\label{rem:ZDelta-dual}
For each $\beta\in \varDelta_+$ fix a non zero element $\xi_{\beta}\in\n_+$ of degree $\beta$. 
Then $\{\xi_{\beta}|\beta\in\varDelta_+\}$ is a basis of $\n_+$; moreover, if $\beta,\gamma\in\varDelta_{+}$ are such that $\beta+\gamma\in\varDelta_+$, then there exists $\coefd{\beta}{\gamma}\ne 0$ such that $[\xi_{\beta},\xi_{\gamma}]=\coefd{\beta}{\gamma}\xi_{\beta+\gamma}$, and if $\beta+\gamma\notin\varDelta_+$, then $[\xi_{\beta},\xi_{\gamma}]=0$.
\end{remark}

Next we introduce the family of subsets of $\varDelta_{+}$ which will parametrize the graded Hopf ideals of $\Zc_{\varDelta}$.

\begin{definition}\label{def:compatible-subsets}
Let $A, B \subseteq \varDelta_{+}$.
\begin{enumerate}[leftmargin=*,label=\rm{(\roman*)}]
\item We say that $A$ is \emph{compatible} if for all $\gamma\in A$ and all pair $\alpha,\beta\in\varDelta_{+}$ such that $\gamma=\alpha+\beta$, then either  $\alpha\in A$ or $\beta\in A$.
\item We say that $B$ is \emph{closed by sums} if for all $\alpha,\beta\in B$, either $\alpha+\beta\notin\varDelta_+$ or $\alpha+\beta\in B$.
\end{enumerate}
\end{definition}

\begin{remark}\label{rem:compatible-closed-by-sum}
$A$ is compatible if and only if $\varDelta_+-A$ is closed by sums.
\end{remark}

We denote by $\comp{\varDelta_+}$ the set of all subsets of $\varDelta_{+}$ that are closed by sums. 

\begin{remark}
$\comp{\varDelta_+}$ is a subposet of $\Pc(\varDelta_+)$ with maximum $\varDelta_+$ and minimum $\emptyset$, as both trivial subsets are closed by sums.
\end{remark}

To deal with subsets of roots as above we need first a statement on sums of roots.

\begin{lemma}\label{lem:sum-positive-roots}
Let $m\ge 3$.
If $\alpha,\gamma_1,\cdots,\gamma_m\in\varDelta_+$ are such that $\alpha=\sum_{i=1}^m\gamma_i$, then there exist $j<k\in\I_m$ such that $\gamma_j+\gamma_k\in\varDelta_+$.
\end{lemma}
\pf
By induction on $\he(\alpha)\ge 3$. If $\he(\alpha)=3$, then either $\alpha=2\alpha_r+\alpha_s$, $r\ne s$ or else $\alpha=\alpha_r+\alpha_s+\alpha_t$: In both cases $m=3$ and the $\gamma_i$'s are simple roots: in the first case $\alpha_r+\alpha_s\in\varDelta_+$, while in the second $\alpha_r+\alpha_s\in\varDelta_+$ up to permute the sub-indices. 

Now assume that the statement holds for roots of height $\le h$ and take $\alpha,\gamma_1,\cdots,\gamma_m\in\varDelta_+$ such that $\he(\alpha)=h+1$, $\alpha=\sum_{i=1}^m\gamma_i$. As the Cartan matrix of $\varDelta$ is finite, there exists $\ell\in\I$ such that $\alpha_{\ell}^{\vee}(\alpha)>0$, hence
$\he(s_{\ell}(\alpha))<\he(\alpha)$. We have two possibilities:
\begin{itemize}[leftmargin=*]
\item $\gamma_i\ne\alpha_{\ell}$ for all $i\in\I_m$. Thus $s_{\ell}(\gamma_i),s_{\ell}(\alpha)\in\varDelta_+$, and $s_{\ell}(\alpha)=\sum_{i=1}^ms_{\ell}(\gamma_i)$. Applying inductive hypothesis for $s_{\ell}(\alpha)$, there exist $j<k\in\I_m$ such that $s_{\ell}(\gamma_j)+s_{\ell}(\gamma_k)\in\varDelta_+$. As $s_{\ell}(\gamma_j)+s_{\ell}(\gamma_k)=s_{\ell}(\gamma_j+\gamma_k)$ and $\gamma_j+\gamma_k\ne\alpha_{\ell}$ since $\he(\gamma_j+\gamma_k)\ge 2$, we have that $\gamma_j+\gamma_k\in\varDelta_+$.
\item There exists $i$ such that $\gamma_i=\alpha_{\ell}$. Up to permute the indices we may assume that $\gamma_m=\alpha_{\ell}$. 
We know that $\alpha-k\alpha_{\ell}\in\varDelta_+$ for all $0\le k\le \alpha_{\ell}^{\vee}(\alpha)$; in particular, $\alpha-\alpha_{\ell}\in\varDelta_+$ and $\alpha-\alpha_{\ell}=\sum_{i=1}^{m-1}\gamma_i$.
If $m=3$, then $\gamma_1+\gamma_2=\alpha-\alpha_{\ell}\in\varDelta_+$. If $m>3$, then we apply inductive hypothesis for $\alpha-\alpha_{\ell}\in\varDelta_+$.
\end{itemize}
In any case, there exist $j<k\in\I_m$ such that $\gamma_j+\gamma_k\in\varDelta_+$.
\epf

We write a slightly different version of compatibility which will be useful in the forthcoming results.

\begin{prop}\label{prop:compatible-alternative}
A subset $A$ is compatible if and only if for all $\alpha\in A$ and $\gamma_i\in\varDelta_{+}$ such that 
$\alpha=\sum_{i=1}^n\gamma_i$, then there exists $i$ such that $\gamma_i\in A$.
\end{prop}

\pf
$(\Leftarrow)$ The case $n=2$ is exactly the definition of compatibility.

\noindent
$(\Rightarrow)$ By induction on $n$: If $n=2$, then it holds by definition. If the statement holds for sums of less than $n$ positive roots and $\alpha=\sum_{i=1}^n\gamma_i$, $n\ge 3$, with $\gamma_i\in\varDelta_{+}$ then we can apply Lemma \ref{lem:sum-positive-roots}: up to permutation we may assume that $\gamma':=\gamma_{n-1}+\gamma_n\in\varDelta_+$. But $\alpha=\sum_{i=1}^{n-2}\gamma_i+\gamma'$ and we can apply inductive hypothesis: either $\gamma_i\in A$ for some $i\le n-2$ (in which case we are done) or else $\gamma'\in A$, which says that either $\gamma_{n-1}\in A$ or $\gamma_n\in A$ by definition of compatibility.
\epf

Now we check that subsets closed by sums classify Lie subalgebras of $\n_+$. More explicitly:

\begin{prop}\label{prop:closed-sums-Lie-subalgs}
There exists a bijective correspondence between subsets closed by sums and $Q_+$-graded Lie subalgebras of $\n_+$ given by
$B\longmapsto \n(B):=\oplus_{\beta\in B} \Bbbk \xi_{\beta}$.
\end{prop}
\pf
First we check that the map is well defined, i.e. each $\n(B)$ is a Lie subalgebra. This follows by Remark \ref{rem:ZDelta-dual}. Indeed, for each pair $\beta,\gamma\in B$, either $\beta+\gamma\notin\varDelta_{+}$, in which case $[\xi_{\beta},\xi_{\gamma}]=0$, or else $\beta+\gamma\in B$ and $[\xi_{\beta},\xi_{\gamma}]=\coefd{\beta}{\gamma}\xi_{\beta+\gamma}\in\n(B)$.

The map is injective by definition, thus it remains to check that it is surjective. Let $\n\subseteq\n_+$ be a $Q_+$-graded Lie subalgebra. As $\n$ is $Q_+$-graded, $\n=\oplus_{\beta\in B}\Bbbk\xi_{\beta}$, where $B$ is the subset of non-trivial homogeneous components. We have to check that $B$ is closed by sums, which follows again by Remark \ref{rem:ZDelta-dual} as $\n$ is a Lie subalgebra.
\epf

Next we introduce a family of ideals and quotients of $\Zc_{\varDelta}$ indexed by $\Pc(\varDelta_+)$. For each
$B \subseteq \varDelta_+$, then we set 
\begin{align}\label{eq:I(A)}
I(B)& :=\langle x_{\beta}|\beta\in \varDelta_+-B\rangle, & Z(B)&:=\Zc_{\varDelta}/I(B).
\end{align}
By definition, $Z(B)$ is a polynomial ring in variables (the images of) $z_{\beta}$, $\beta\in B$.
\smallbreak

We are mostly interested in those ideals attached to $B\in\comp{\varDelta_{+}}$ as we will see next.

\begin{lemma}\label{lem:good-definition-of-I}
If $B\in\comp{\varDelta_{+}}$, then $I(B)$ is a $Q_+$-graded Hopf ideal. 
\end{lemma}
\pf
Let $A=\varDelta_+-B$. It is enough to prove that $\Delta(z_{\alpha})\in I(B)\ot \Zc_{\varDelta}+\Zc_{\varDelta}\ot I(B)$ for all $\alpha\in A$. Let $\alpha\in A$:
$\Delta(z_{\alpha})-z_{\alpha}\ot 1-1\otimes z_{\alpha}$ is a linear combination of terms 
\begin{align*}
&z_{\gamma_1}\cdots z_{\gamma_k}\ot z_{\gamma_{k+1}}\cdots z_{\gamma_m}, &&1\le k<m, \, \gamma_i\in\varDelta_{+}, \, \sum_{i=1}^{m}\gamma_i=\alpha.
\end{align*}
By Proposition \ref{prop:compatible-alternative}, for each term $z_{\gamma_1}\cdots z_{\gamma_k}\ot z_{\gamma_{k+1}}\cdots z_{\gamma_m}$ there exists $i$ such that $\gamma_i\in A$, thus it belongs to $I(B)\ot \Zc_{\varDelta}+\Zc_{\varDelta}\ot I(B)$.
\epf

We will use the previous result to give a parametrization of graded Hopf ideals of $\Zc_{\varDelta}$.

\begin{theorem}\label{thm:Z-Delta-poset-ideals-sets}
There exists an anti-isomorphism of posets between $\comp{\varDelta_+}$ and the set of $Q_+$-graded Hopf ideals of $\Zc_{\varDelta}$ given by
$$ B\in \comp{\varDelta_+} \longmapsto Z(B). $$
\end{theorem}
\pf
By Lemma \ref{lem:good-definition-of-I} the map is well defined. Moreover, the map is an anti-morphism of posets, and is injective since if $B \neq B'\in \comp{\varDelta_+}$, then we may assume that there exists $\alpha\in B'-B$, in which case $z_{\alpha}\in I(B)-I(B')$.

Thus, it remains to prove that the map is surjective. For this, we first notice the following. For each $B\in\comp{\varDelta_+}$, as $I(B)$ is a Hopf ideal, $I(B)^{\perp}\subseteq U(\n_+)$\footnote{As before, here we take the $Q^+$-graded dual of $\Zc_{\varDelta}$.} is a Hopf subalgebra of $U(\n_+)$ by \cite[Proposition 5.2.5]{Rad-libro}. Thus $I(B)^{\perp}=U(\n)$ for some Lie subalgebra $\n$ of $\n_+$. By Proposition \ref{prop:closed-sums-Lie-subalgs}, $\n=\n(B')$ for some subset $B'\subseteq \varDelta_+$ closed by sums: clearly, $B=B'$. 
Now take $I$ a $Q_+$-graded Hopf ideal of $\Zc_{\varDelta}$: again, the subspace $I^{\perp}\subseteq U(\n_+)$ is a Hopf subalgebra of $U(\n_+)$, so $I^{\perp}=U(\n(C))$ for some subset $C\subseteq \varDelta_+$ closed by sums. But by the previous argument, $U(\n(C))=I(B)^{\perp}$ for $B=C$, hence $I=I(B)$.
\epf

\section{Posets of pre-Nichols algebras of diagonal type}\label{sec:poset-pre-Nichols}

Here we proceed to describe the poset of all graded pre-Nichols algebras of $\bq$ with finite $\GK$. 
First we assume that the diagram of $\bq$ is connected and consider two cases: whether the eminent pre-Nichols algebra is or is not the distinguished one.
Later on we give an approach towards the non-connected case: we have an obstruction to attack the general case coming from those cases where the eminent pre-Nichols algebra is not known.

\subsection{Eminent pre-Nichols algebras which are not distinguished}\label{subsec:eminent-non-dist}

Throughout this section, $\bq$ will denote a braiding matrix of one of the following types: Cartan $A_2$ with $q\in \G_3'$,
$\superqa{3}{q}{\{2\}}$ with $q\in \G_{\infty}$, $\superqa{3}{q}{\{1,2,3\}}$ with $q\in \G_{\infty}$,
$\g(2,3)$ with diagram $d_1$ or $\g(2,3)$ with diagram $d_2$. Therefore, the eminent pre-Nichols algebra $\htoba_{\bq}$ is not the distinguished pre-Nichols algebra $\wtoba_{\bq}$ as described in Theorem \ref{thm:eminent-not-distinguished}, and we can consider three subalgebras of coinvariants, associated to the non-trivial canonical projections
$\htoba_{\bq}\twoheadrightarrow \wtoba_{\bq} \twoheadrightarrow \toba_{\bq}$:
\begin{align}\label{eq:coinv-subalgebras}
\hZ_{\bq}&=\htoba_{\bq}^{\co \toba_{\bq}},  & \Zc_{\bq}&=\wtoba_{\bq}^{\co \toba_{\bq}}, & 
Z_{\bq}&=\htoba_{\bq}^{\co \wtoba_{\bq}}
\end{align}
By Theorem \ref{thm:eminent-not-distinguished}, $Z_{\bq}$ is a skew central Hopf subalgebra of $\htoba_{\bq}$, and by \cite{An-dist} $\Zc_{\bq}$ is a skew central Hopf subalgebra of $\wtoba_{\bq}$, and both are skew polynomial algebras: $\underline{\fO^{\bq}}$ is the set of degrees of the generators of $\Zc_{\bq}$. 
We will see that $\hZ_{\bq}$ is also polynomial algebras whose generators are obtained by joining the generators of $\Zc_{\bq}$ and 
$Z_{\bq}$, and at the same time a skew central Hopf subalgebra of $\htoba_{\bq}$. To this end we introduce the set $\hOc{\bq}$ by extending $\underline{\fO^{\bq}_+}$ with the degrees of the generators of $Z_{\bq}$. Explicitly,
\begin{align}\label{eq:extended-root-system}
\hOc{\bq} &:=\begin{cases}
\{1^3, 1^32^3,2^3, \mathbf{1^22},\mathbf{12^2}\}, & \text{type }A_2, q\in \G_3',
\\
\{1^N, 3^N, \mathbf{12^23}\}, & \text{type }\superqa{3}{q}{\{2\}}, q\in \G_{\infty},
\\
\{1^N2^N, 2^N3^N, \mathbf{13}\}, & \text{type }\superqa{3}{q}{\{1,2,3\}}, q\in \G_{\infty},
\\
\{1^32^3,1^32^63^3,2^33^3,1^62^63^6, \mathbf{1^22^33}, \mathbf{12^33^2}\}, & \text{type }\g(2,3) \text{ with diagram }d_1,
\\
\{1^32^33^3,1^32^63^3,2^3,2^63^6, \mathbf{12^33}, \mathbf{1^22^33^2}\}, & \text{type }\g(2,3) \text{ with diagram }d_2,
\end{cases}
\end{align}
where bold degrees are those of the generators of $Z_{\bq}$. By \cite[Lemma 2.4]{ACSa},
\begin{align*}
&
\left.
\begin{aligned}
\Hc_{\htoba_{\bq}} &= \Hc_{Z_{\bq}} \Hc_{\wtoba_{\bq}}= \Hc_{\hZ_{\bq}} \Hc_{\toba_{\bq}} \\ 
\Hc_{\wtoba_{\bq}} &= \Hc_{\Zc_{\bq}} \Hc_{\toba_{\bq}}
\end{aligned}
\right\}
&& \implies & \Hc_{\hZ_{\bq}} &=\Hc_{Z_{\bq}}\Hc_{\Zc_{\bq}}=\prod_{\beta\in \hOc{\bq}} \tfrac{1}{1-t^{\beta}}.
\end{align*}

\begin{prop}\label{prop:eminent-non-dist-poset}
\begin{enumerate}[leftmargin=*,label=\rm{(\roman*)}]
\item\label{item:eminent-non-dist-central-Hopf} $\hZ_{\bq}$ is a skew central Hopf subalgebra of $\htoba_{\bq}$, and is a skew polynomials ring whose generators are homogeneous: $\hOc{\bq}$ is the set of their degrees.
\item\label{item:eminent-non-dist-central-pre-Nichols} For each subset $B\subseteq \hOc{\bq}$ closed by sums, there exists a $\zt$-graded pre-Nichols algebra $\toba(\bq,B)$ with Hilbert series
\begin{align*}
\Hc_{\toba(\bq,B)}(t)&=\Hc_{\toba_{\bq}}(t) \left(\prod_{\alpha\in B}\tfrac{1}{1-t^{\alpha}}\right).
\end{align*}
\item\label{item:eminent-non-dist-central-poset-isom} The map $\comp{\hOc{\bq}}\to \prefdgr{\bq}$, $B\mapsto \toba(\bq,B)$ is an anti-isomorphism of posets.
\end{enumerate}

\end{prop}
\pf
We prove the three statements for each case. First assume that $\bq$ is of Cartan type  $A_2$ with $q\in \G_3'$: by \cite[Lemma 4.10]{ASa}, $\hZ_{\bq}$ is a central Hopf subalgebra and, as an algebra, is a skew polynomial ring in variables $x_1^3$, $x_2^3$, $x_{112}$, $x_{221}$ and $x_{12}^3$, so \ref{item:eminent-non-dist-central-Hopf} holds. For \ref{item:eminent-non-dist-central-pre-Nichols}, notice that $\alpha,\beta \in \hOc{\bq}$ are such that $\alpha+\beta\in \hOc{\bq}$ iff $\alpha+\beta=1^32^3$ and either $\{\alpha,\beta\}=\{1^3,2^3\}$ or else $\{\alpha,\beta\}=\{1^22,12^2\}$.
By \cite[Lemma 4.10]{ASa}, $x_1^3$, $x_2^3$, $x_{112}$, $x_{221}$ are primitive and the following formula holds in $\htoba_{\bq}$:
\begin{align}\label{eq:coproduct-A2-Cartan}
\Delta(x_{12}^3) &= x_{12}^3\ot 1 + 1\ot x_{12}^3+(1-q^2)q_{21}^3 x_1^3 \otimes x_2^3 + (q^2-q)x_{112}\ot x_{221}.
\end{align}
Thus, if $B\subseteq \hOc{\bq}$ is closed by sums, then $\toba(\bq,B)=\htoba_{\bq}/\langle x_{\alpha}, \alpha \in \hOc{\bq}-B \rangle$ is a pre-Nichols algebra with the desired Hilbert series. Reciprocally, if $\toba$ is a $\zt$-graded pre-Nichols algebra, then either $\toba=\htoba_{\bq}$ or else one of the primitive elements $x_1^3$, $x_2^3$, $x_{112}$, $x_{221}$ is zero, since the subspace of primitive elements is spanned by these primitive elements (not in degree one) and $V$. Let $T_1$ be the set of degrees of $\{x_1^3, x_2^3, x_{112}, x_{221}\}$ of those elements annihilating in $\toba$. If $T_1\cap \{x_1^3, x_2^3\}=\emptyset$ or $T_1\cap \{x_{112}, x_{221}\}=\emptyset$ then $x_{12}^3$ cannot be zero in $\toba$, otherwise $x_{12}^3$ is primitive and this element may be zero in $\toba$. Thus set $T=T_1$ if $x_{12}^3\ne 0$ in $\toba$, or $T=T_1\cup\{1^32^3\}$ if $x_{12}^3=0$ in $\toba$. Then $B:=T^c$ is closed by sums, with $\toba=\toba(\bq,B)$, and hence \ref{item:eminent-non-dist-central-poset-isom} follows.

\medbreak
Next assume that $\bq$ is of type $\superqa{3}{q}{\{2\}}$ with $q\in \G_{\infty}$. By the proof of \cite[Proposition 5.5]{ACSa}, $x_{12^23}$ is skew central. Using the relations verified in that proof, also $x_1^N$ and $x_3^N$ are skew central, and the three elements are primitive. Thus $\hZ_{\bq}$ is a skew polynomial algebra with generators $x_{12^23}$, $x_1^N$ and $x_3^N$, and \ref{item:eminent-non-dist-central-Hopf} follows. Now every subset of $\hOc{\bq}$ is closed by sums since $\alpha+\beta\notin \hOc{\bq}$ for all $\alpha,\beta \in \hOc{\bq}$. Thus the proof of \ref{item:eminent-non-dist-central-pre-Nichols} and \ref{item:eminent-non-dist-central-poset-isom} are straightforward. The case in which $\bq$ is of type $\superqa{3}{q}{\{1,2,3\}}$ with $q\in \G_{\infty}$ follows analogously taking into account the proof of \cite[Proposition 5.6]{ACSa}.

\medbreak
If $\bq$ is of type $\g(2,3)$ with diagram $d_1$, then
\begin{align*}
x_u &=[[x_{12},x_{123}]_c, x_2]_c && \text{and} & x_v&=[[x_{123},x_{23}]_c, x_2]_c
\end{align*}
are primitive in $\htoba_{\bq}$ and skew central by \cite[Proposition 4.2]{ACSa2}. Also,
\begin{align*}
	[x_i,x_{23}^3]_c&=[x_i,x_{12}^3]_c=0 &&\text{ for all }i\in\{1,2,3\}
\end{align*}
Indeed, the proof of that proposition says that $[x_1,x_{23}^3]_c=[x_3,x_{12}^3]_c=0$, and the other relations follow from the quantum Serre relations. Thus $x_{12}^3$, $x_{23}^3$ are skew-central in $\htoba_{\bq}$. As these elements are primitive in $\wtoba_{\bq}$, 
\begin{align*}
\Delta(x_{ij}^3) &\in x_{ij}^3 \ot 1 + 1 \ot x_{ij}^3 + \langle x_u, x_v \rangle \ot \htoba_{\bq} + \htoba_{\bq} \ot \langle x_u, x_v \rangle.
\end{align*}
Taking into account the $\Z^{\theta}$-degree we check that $x_{12}^3$ and $x_{23}^3$ are primitive in $\htoba_{\bq}$ as well.

Now we claim that $B=\{x_u^{a_1} x_v^{a_2}  x_{12}^{3a_3}x_{12^23}^{3a_4}x_{23}^{3a_5}x_{123}^{6a_6} | a_i\in\N_0\}$ is a basis of $\hZ_{\bq}$. Indeed these elements belong to $\hZ_{\bq}$, are linearly independent by \cite[Proposition 4.2]{ACSa2} and then they must generate $\hZ_{\bq}$ because of the expression of $\Hc_{\hZ_{\bq}}$ above. As $\hZ_{\bq}$ is a normal subalgebra of $\htoba_{\bq}$, we have that $\ad_c x_i(x_{12^23}^3)\in\hZ_{\bq}$ for $i=1,2,3$. But we can check that there are no elements of degrees $1^42^63^3$, $1^32^73^3$, $1^32^63^4$ in $B$, so $\ad_c x_i(x_{12^23}^3)=0$ for all $i=1,2,3$. A similar argument shows that $\ad_c x_i(x_{123}^6)=0$ for all $i=1,2,3$, so $x_{12^23}^3$ and $x_{123}^6$ are skew central. It means that $\hZ_{\bq}$ is skew central, and also a coideal subalgebra on both sides, so it is a central Hopf subalgebra.
As $\hZ_{\bq}$ is a skew polynomial algebra with generators $x_{12}^3$, $x_{12^23}^3$, $x_{23}^3$, $x_{123}^6$, $x_u$ and $x_v$,  \ref{item:eminent-non-dist-central-Hopf} follows.

Now $x_{123}^6$ is skew primitive since the unique pairs of elements $(b_1,b_2)\in B\times B$ such that the sums of their degrees is $1^62^63^6$ are $(x_{123}^6,1)$ and $(1,x_{123}^6)$.  Now we compute $\Delta(x_{12^23}^3)$. By direct computation,
\begin{align*}
\Delta(x_{12^23}) &= x_{12^23} \ot 1 + 1 \otimes x_{12^23} + 3\zeta q_{23} x_{12}\ot x_{23} + (1-\zeta^2) x_{123}\ot x_2+ 3 x_1\ot x_{23}x_2.
\end{align*}
Using the degree again, the possible non-trivial summands of $\Delta(x_{12^23}^3)$ are $x_{12}^3\ot x_{23}^3$ and $x_u \ot x_v$.  Hence we just look at the corresponding degrees and use that $x_u=[x_{12},x_{12^23}]_c$ and $x_v=[x_{12^23},x_{23}]_c$ to check the following identity:
\begin{align}\label{eq:coproduct-g23-d1}
\Delta(x_{12^23}^3) &= x_{12^23}^3 \ot 1 + 1 \otimes x_{12^23}^3 -27 q_{21}^3q_{31}^3 x_{12}^3\ot x_{23}^3 + 3\zeta q_{21}^2q_{31}^2q_{32} x_u \ot x_v.
\end{align}
Next we see that $A \subseteq \hOc{\bq}$ is compatible if and only if either $1^32^63^3\notin A$ or else $1^32^63^3\in A$ and $\{1^32^3,2^33^3\}\cap A, \{1^22^33,12^33^2\}\cap A \neq \emptyset$, and a subset is compatible if and only of the complement is closed by sums, as in Remark \ref{rem:compatible-closed-by-sum}. 
Thus \ref{item:eminent-non-dist-central-pre-Nichols} and \ref{item:eminent-non-dist-central-poset-isom} follow by an argument analogous to the Cartan case $A_2$.

\medbreak
Finally, assume that $\bq$ is of type $\g(2,3)$ with diagram $d_2$. The proof of this case is analogous to the one above but using \cite[Proposition 4.3]{ACSa2} and its proof. Indeed, set
\begin{align*}
x_u &=[x_{12^23},x_2]_c && \text{and} & x_v&=[x_{123},x_{12^23}]_c.
\end{align*}
By loc. cit. $x_u$ and $x_v$ are primitive in $\htoba_{\bq}$ and skew central. In addition, $x_2^3$ is primitive and skew central by direct computation, and 
$x_{23}^6$ and $x_{123}^3$ are also primitive and skew-central in $\htoba_{\bq}$ since they are primitive and skew central in $\wtoba_{\bq}$ and we take into account the $\Z^{\theta}$-degree as above. Next we observe that $B=\{x_u^{a_1} x_v^{a_2}  x_{123}^{3a_3}x_{12^23}^{3a_4}x_{2}^{3a_5}x_{23}^{6a_6} | a_i\in\N_0\}$ is a basis of $\hZ_{\bq}$ by looking at the expression of the Hilbert series $\Hc_{\hZ_{\bq}}$ as above, so $\hZ_{\bq}$ is a central Hopf subalgebra. As $\hZ_{\bq}$ is a skew polynomial algebra with generators $x_{2}^3$, $x_{12^23}^3$, $x_{123}^3$, $x_{23}^6$, $x_u$ and $x_v$, \ref{item:eminent-non-dist-central-Hopf} follows. 

We need now an explicit expression of $\Delta(x_{12^23}^3)$. To do so we compute first
\begin{align*}
\Delta(x_{12^23}) &= x_{12^23} \ot 1 + 1 \otimes x_{12^23} -q_{23}(1-\zeta^2) x_{12}\ot x_{23} 
\\
& \qquad + (1-\zeta^2) x_{123}\ot x_2-3\zeta^2 x_1\ot x_{23}x_2.
\end{align*}
Working as in \eqref{eq:coproduct-g23-d1} we obtain the following:
\begin{align}\label{eq:coproduct-g23-d2}
\Delta(x_{12^23}^3) &= x_{12^23}^3 \ot 1 + 1 \otimes x_{12^23}^3 + 3q_{21}^3q_{23}^3\zeta(1-\zeta) x_{123}^3\ot x_{2}^3 + q_{21}^2q_{23}^2(\zeta-1) x_v \ot x_u.
\end{align}
Now, $A \subseteq \hOc{\bq}$ is compatible if and only if either $1^32^63^3\notin A$ or else $1^32^63^3\in A$ and $\{1^32^33^3,2^3\}\cap A, \{1^22^33^2,12^33\}\cap A \neq \emptyset$, so \ref{item:eminent-non-dist-central-pre-Nichols} and \ref{item:eminent-non-dist-central-poset-isom} follow as above.
\epf

We finish this subsection by identifying the algebraic groups of types Cartan $A_2$ and $\g(2,3)$ (for both diagrams). Assume that $\bq$ is such that $\hZ_{\bq}$ is a central Hopf subalgebra (that we can assume up to twist the braiding by Lemma \ref{lem:twist-equivalent-central}).

Let $\n_+$ be the positive part of the Lie algebra $\g=\sli_5$, which is a 10-dimensional nilpotent algebra with generators $e_{ij}$, $1\le i\le j\le 4$, where
\begin{align*}
e_{ii}&=e_i, & e_{ij}&=[e_{ik},e_{k+1\, j}], &&i\le k<j.
\end{align*}
Notice that $e_{ij}$ has degree $\alpha_{ij}$ for all $i\le j$. 

Let $Z_{\varDelta}$ be the corresponding algebraic group with Lie algebra $\n_+$.  The subset
\begin{align*}
 B=\{1, 123, 1234,234,4\} \subseteq \varDelta_{+}
\end{align*}
is closed by sums. The associated quotient Hopf algebra $Z(B)$ in Theorem \ref{thm:Z-Delta-poset-ideals-sets} is a polynomial ring in variables $z_{\beta}$, $\beta\in B$, where all $z_{\beta}\ne z_{1234}$ are primitive, and
\begin{align*}
	\Delta(z_{1234}) &= z_{1234}\ot 1+ 1\ot z_{1234}+z_{1}\ot z_{234}+z_{123}\ot z_{4}.
\end{align*}

\begin{lemma}\label{lem:identifying-hZ-new} 
Assume that $\bq$ is such that $\hZ_{\bq}$ is a central Hopf subalgebra.
\begin{enumerate}[leftmargin=*,label=\rm{(\roman*)}]
\item If $\bq$ is of Cartan type $A_2$ and $q\in\G_3'$, then $\hZ_{\bq}\simeq Z(B)$ as Hopf algebras.
\item If $\bq$ is of type $\superqa{3}{q}{\{2\}}$ or $\superqa{3}{q}{\{1,2,3\}}$ with $q\in \G_{\infty}$, then $\hZ_{\bq}\simeq \Bbbk[z_1,z_2,z_3]$ as Hopf algebras.
\item If $\bq$ is of type $\g(2,3)$ and diagram $d_1$ or $d_2$, then $\hZ_{\bq}\simeq Z(B)\times \Bbbk[z]$ as Hopf algebras.
\end{enumerate}
\end{lemma}
\pf
It follows case-by-case, using the coproduct formulas \eqref{eq:coproduct-A2-Cartan}, \eqref{eq:coproduct-g23-d1} and \eqref{eq:coproduct-g23-d2} in the proof of Proposition \ref{prop:eminent-non-dist-poset}.
\epf

\subsection{The connected case}\label{subsec:poset}

Let $\bq$ be a braiding of diagonal type such that $\dim \toba_{\bq}<\infty$, the Dynkin diagram is connected 
and $\bq$ is not of Cartan type $A_{\theta}$, $D_{\theta}$ with $q=-1$.

Those cases where the eminent pre-Nichols algebra $\htoba_{\bq}$ is not the distinguished pre-Nichols algebra $\wtoba_{\bq}$ were treated in \S \ref{subsec:eminent-non-dist}, so assume for a while that $\htoba_{\bq}=\wtoba_{\bq}$. Thus $\hZ_{\bq}=\htoba_{\bq}^{\co \toba_{\bq}}=\Zc_{\bq}$ is a skew central Hopf subalgebra of $\htoba_{\bq}$, and an skew polynomial algebra such that $\underline{\fO^{\bq}_+}$ is the set of degrees of the generators of $\Zc_{\bq}$. Hence we set $\hOc{\bq}=\underline{\fO^{\bq}_+}$ for this case: we have already defined $\hOc{\bq}$ for the other five kinds of braidings in \eqref{eq:extended-root-system}, so in any case $\hOc{\bq}$ is the set of degrees of the generators of $\hZ_{\bq}$.

We will extend Proposition \ref{prop:eminent-non-dist-poset} to any $\bq$ with connected Dynkin diagram. We deal first with the existence of pre-Nichols algebras. For each $\ubeta\in\hOc{\bq}$ let $z_{\ubeta}$ be the corresponding generator of degree $\ubeta$: if $\ubeta=N_{\beta}\beta$ for some Cartan root $\beta$, then $z_{\ubeta}:=x_{\beta}^{N_{\beta}}$, otherwise $z_{\ubeta}$ is the extra relation of such degree. For example, for Cartan type $A_2$ and $q\in\G_3'$, $z_{1^22}=x_{112}$ and $z_{12^2}=x_{221}$.

\begin{lemma}\label{lemma:compatible-sets-to-quotients}\label{lem:subset-closed-sums-to-pre-Nichols}
Let $\bq$ be a matrix with connected Dynkin diagram such that $\dim\toba_{\bq}<\infty$ and is not of Cartan type $A_n$, $D_n$ with $q=-1$. 

For each $B\subseteq \hOc{\bq}$ closed by sums, the quotient
\begin{align*}
\toba(\bq,B) := \htoba_{\bq}/\langle z_{\ubeta} | \ubeta\in \hOc{\bq}-B \rangle
\end{align*}
is a $\zt$-graded pre-Nichols algebra $\toba(\bq,B)$ with Hilbert series
\begin{align}\label{eq:compatible-sets-to-quotients-Hilbert-series}
\Hc_{\toba(\bq,B)}(t)&=\Hc_{\toba_{\bq}}(t) \left(\prod_{\ubeta\in B}\tfrac{1}{1-t^{\ubeta}}\right).
\end{align}
\end{lemma}
\pf
As mentioned above, it suffices to deal with the case $\htoba_{\bq}=\wtoba_{\bq}$ as in Theorem \ref{thm:eminent-distinguished}, since the statement holds for other $\bq$ by Proposition \ref{prop:eminent-non-dist-poset} \ref{item:eminent-non-dist-central-pre-Nichols}. Here, $\hOc{\bq}=\varDelta_+$ is the set of positive roots of the semisimple Lie algebra attached to $\wtoba_{\bq}$ in \cite{AAR-imrn}. Moreover, we can assume that $\hZ_{\bq}$ is central up to change $\bq$ by a twist equivalent braiding matrix, see Proposition \ref{prop:twist-equiv-poset-preNichols} and Lemma \ref{lem:twist-equivalent-central}, in which case $\hZ_{\bq}=\Zc_{\varDelta}$.

Let $B \subseteq \hOc{\bq}$ be a subset closed by sums. By Remark \ref{rem:compatible-closed-by-sum}, $A=B^c$ is compatible, so $I(B)=\langle z_{\ubeta} |\ubeta \in A \rangle$ is an $\N_0^\theta$-graded Hopf ideal of $\Zc_{\varDelta}$ by Lemma \ref{lem:good-definition-of-I}. Thus $\htoba_{\bq}I(B)$ is a Hopf ideal of $\htoba_{\bq}$. By Proposition \ref{prop:HopfidealB-HopfidealZ} \ref{item:HopfidealB-HopfidealZ-graded}, $\Hc_{\wtoba_{\bq}}(t)=\Hc_{\toba(\bq,B)}(t)\Hc_{\mathtt{Z}}(t)$, where $\mathtt{Z}$ is the subalgebra of coinvariants of the projection $\Zc_{\bq}\twoheadrightarrow \Zc_{\bq}/I(B)$. As the subalgebra $Z$ generated by $z_{\ubeta}$, $\ubeta\in A$, is simultaneously a coideal subalgebra and  polynomial ring in these variables such that $\langle Z^+\rangle=I(B)$, we have that $Z=\mathtt{Z}$  by Proposition \ref{prop:HopfIdeal-Subcoalgebras}, so 
$\Hc_{\mathtt{Z}}(t)=\left(\prod_{\alpha\in A}\tfrac{1}{1-t^{\alpha}}\right)$. As
\begin{align*}
\Hc_{\wtoba_{\bq}}(t)=\Hc_{\toba_{\bq}}(t)\left(\prod_{\alpha\in \hOc{\bq}}\tfrac{1}{1-t^{\alpha}}\right)
\end{align*}
by \eqref{eq:Hilbert-series-distinguished}, we obtain the expression \eqref{eq:compatible-sets-to-quotients-Hilbert-series} for $\Hc_{\toba(\bq,B)}(t)$.
\epf

\begin{remark}\label{rem:PBW-basis-pre-Nichols}
Let $B=\{\ugamma_1,\cdots,\ugamma_{L}\}$ be a numeration of $B$. Then the set
\begin{align*}
&x_{\beta_1}^{n_1} \cdots x_{\beta_M}^{n_M} z_{\ugamma_1}^{p_1}\cdots z_{\ugamma_L}^{p_L}, && 0\le n_i <N_{\beta_i}, \, 0 \le p_j<\infty,
\end{align*}
is a basis of $\toba(\bq,B)$.
\end{remark}

Now we state a characterization of the poset of $\N_{0}^{\theta}$-graded pre-Nichols algebras.

\begin{theorem}\label{thm:connected-isom-posets}
Let $\bq$ a matrix with connected Dynkin diagram such that $\dim\toba_{\bq}<\infty$ and is not of Cartan type $A_n$, $D_n$ with $q=-1$. The map 
\begin{align*}
& \comp{\hOc{\bq}} \to \prefdgr{\bq}, && B\mapsto \toba(\bq,B),
\end{align*}
is an anti-isomorphism of posets.
\end{theorem}
\pf
As in the proof of Lemma \ref{lemma:compatible-sets-to-quotients}, we can restrict to those $\bq$ such that $\htoba_{\bq}=\wtoba_{\bq}$ since the remaining cases were treated in Proposition \ref{prop:eminent-non-dist-poset} \ref{item:eminent-non-dist-central-poset-isom}, and we can use Proposition \ref{prop:twist-equiv-poset-preNichols} and Lemma \ref{lem:twist-equivalent-central} to assume that $\Zc_{\bq}$ is central up to change $\bq$ by a twist equivalent matrix.

By Lemma \ref{lemma:compatible-sets-to-quotients}, the map above is injective: if $B\ne B'$ are two different sets closed by sums, then 
$\Hc_{\toba(\bq,B)}(t)\ne \Hc_{\toba(\bq,B')}(t)$. Also, it is a anti morphism of posets.

On the other hand, we will see that the map is also surjective. Let $\toba\in \prefdgr{\bq}$: 
\begin{itemize}[leftmargin=*]
\item By definition, there exists an $\N_0^{\theta}$-graded Hopf ideal $\Ic$ such that $\toba \simeq \htoba_{\bq}/\Ic$.
\item By Proposition \ref{prop:HopfidealB-HopfidealZ} \ref{item:HopfidealB-HopfidealZ-bijection} there exists a graded Hopf ideal $I$ of $\hZ_{\bq}$ such that $\Ic=\htoba_{\bq} I$. 
\item By Theorem \ref{thm:Z-Delta-poset-ideals-sets} there exists a set $B\subseteq \hOc{\bq}=\varDelta_+$ closed by sums such that $I=I(B)$.
\end{itemize}
All in all, we have that 
\begin{align*}
\toba&=\htoba_{\bq}/\htoba_{\bq} I(B)= \htoba_{\bq}/\langle z_{\ubeta} | \ubeta\in A \rangle = \toba(\bq,B),
\end{align*}
so the map is also surjective, and the statement follows.
\epf

\begin{example}
Fix $\zeta\in\G_3'$. Let $\bq$ be a braiding matrix of type $\g(2,3)$ with Dynkin diagram 
\begin{align*}
&d_3: \xymatrix@C=30pt{\overset{-1}{\circ} \ar@{-}[r]^{\xi^2} &\overset{-\xi^2}{\circ}\ar@{-}[r]^{\xi^2}\ &\overset{-1}{\circ}}.
\end{align*}
Then $\wtoba_{\bq}=\htoba$ is presented by generators and relations
\begin{align*}
& x_{13} , & &x_{2221}, & &x_{2223}, & & x_{1}^2, & & x_{3}^2, =0, 
& & [x_1,x_{223}]_c+q_{23} x_{12^23}-(1-\zeta^2)x_2x_{123}.
\end{align*}
Also, $\hOc{\bq}=\{1^22^3, 1^32^63^3, 2^33^3, 2^6\}$ is (isomorphic to) the set of positive roots associated to a Lie algebra of type $A_2\times A_1$, and $\hZ_{\bq}$ is the subalgebra generated by
\begin{align*}
z_1&=x_{12}^3, & z_{12}&=x_{12^23}^3, & z_2&=x_{23}^3, & z_3&=x_{2}^6. 
\end{align*}
The poset $\prefdgr{\bq}$ is the following:
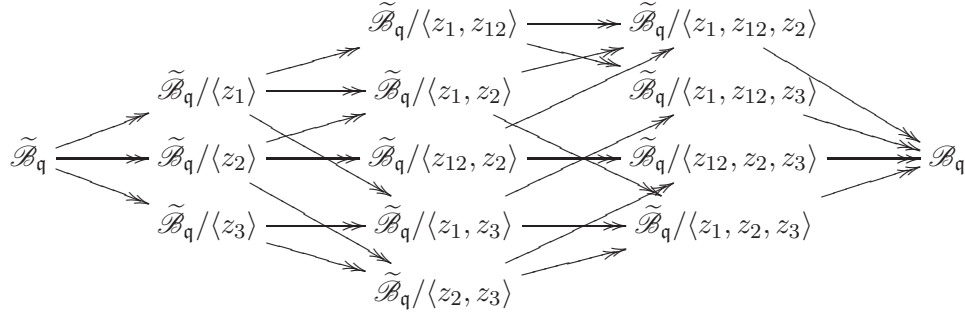
\begin{align*}
\xymatrix@R=5pt@C=35pt{
& & \wtoba_{\bq}/\langle z_1, z_{12} \rangle \ar@{->>}[r] \ar@{->>}[rd] &  \wtoba_{\bq}/\langle z_1, z_{12}, z_2 \rangle \ar@{->>}[rdd]&
\\
& \wtoba_{\bq}/\langle z_1 \rangle \ar@{->>}[r] \ar@{->>}[ru] \ar@{->>}[rdd] & \wtoba_{\bq}/\langle z_1, z_2 \rangle \ar@{->>}[ru] \ar@{->>}[rdd] &  \wtoba_{\bq}/\langle z_1, z_{12}, z_3 \rangle \ar@{->>}[rd]&
\\
\wtoba_{\bq} \ar@{->>}[r] \ar@{->>}[ru] \ar@{->>}[rd] & \wtoba_{\bq}/\langle z_2 \rangle \ar@{->>}[r] \ar@{->>}[ru] \ar@{->>}[rdd] & 
\wtoba_{\bq}/\langle z_{12}, z_2 \rangle \ar@{->>}[ruu] \ar@{->>}[r] &  \wtoba_{\bq}/\langle z_{12}, z_2, z_3 \rangle \ar@{->>}[r]& \toba_{\bq}
\\
& \wtoba_{\bq}/\langle z_3 \rangle \ar@{->>}[r] \ar@{->>}[rd] & \wtoba_{\bq}/\langle z_1, z_3 \rangle \ar@{->>}[ruu] \ar@{->>}[r] &  \wtoba_{\bq}/\langle z_1, z_2, z_3 \rangle \ar@{->>}[ru]&
\\
& & \wtoba_{\bq}/\langle z_2, z_3 \rangle \ar@{->>}[ruu] \ar@{->>}[ru] &  &
}
\end{align*}
When we move from left to right, the $\GK$ goes down from 4 to 0.

We can observe here that the poset of graded pre-Nichols algebras is not preserved by Weyl equivalence: if $\bq$ is also of type $\g(2,3)$ but with Dynkin diagram $d_1$ or $d_2$, then $\prefdgr{\bq}$ has 50 elements, since this is the number of subsets of $\hOc{\bq}$ closed by sums.
\end{example}

\subsection{The non-connected case}

Here we take an arbitraty matrix $\bq$ (whose diagram is not necessarily connected). Following the spirit of \cite[\S 3]{ASa} we set
\begin{align*}
\I^{\pm} &:= \{ i\in \I: q_{ii}=\pm 1, \qti_{ij}=1 \text{ for all }j\ne i\}, & \I^{(c)} &:= \I-\left( \I^+\cup\I^-\right).
\end{align*}
that is, $\I^{\pm}$ contains all the connected components of an isolated vertex labelled with $\pm1$ while $\I^{(c)}$ is the union of those points labelled with $q_{ii}\ne \pm1$ and those connected components with at least two vertices. 
Let $C_1, \cdots, C_d$ be the partition of $\I^{(c)}$ in connected components: i.e. $\I^{(c)}=\cup_{\ell=1}^d C_d$, where the diagram of each $\bq^{(\ell)}:=(q_{ij})_{i,j\in C_{\ell}}$ is connected and $\qti_{ij}=1$ if $i\in C_{\ell_1}$, $j\in C_{\ell_2}$, $\ell_1\ne \ell_2$.

Let $\toba$ be an $\N_0$-graded pre-Nichols algebra of $\bq$.
For each $1\le \ell\le d$, let $\toba^{(\ell)}$ be the subalgebra of $\toba$ generated by $x_i$, $i\in C_{\ell}$, which is a pre-Nichols algebra of $\bq^{(\ell)}$. Similarly, let $\toba^{\pm}$ be the subalgebra of $\toba$ generated by $x_i$, $i\in\I^{\pm}$, and $\bq^{(\pm)}:=(q_{ij})_{i,j\in \I^{\pm}}$.

\begin{lemma}\label{lem:diff-conn-components}
Let $1\le \ell\le d$, $i\in C_{\ell}$, $j\in\I-C_{\ell}$ such that $q_{jj}\ne 1$. If $\GK\toba<\infty$, then $x_{ij}=0$ in $\toba$.
\end{lemma}
\pf
By hypothesis, $i$ and $j$ are not connected, so $x_{ij}$ is primitive in $\toba$. Also, $\ord q_{ii}\ge 2$ since $C_{\ell} \subseteq \I^{(c)}$, and also $\ord q_{jj}\ge 2$ since $q_{jj}\ne 1$. If either $\ord q_{ii}>2$ or $\ord q_{jj}>2$, then $x_{ij}=0$ by \cite[Proposition 3.2]{ASa}. Otherwise $q_{ii}=q_{jj}=-1$ and there exists $k\ne i,j$ such that $\qti_{ik}\ne 1=\qti_{jk}$. Suppose that $x_{ij}\ne 0$ in $\toba$. Set $y_1=x_i$, $y_2=x_k$, $y_3=x_{ij}$ and $W:=\Bbbk y_1+\Bbbk y_2+\Bbbk y_3$ is a three-dimensional subspace of $\Pc(\toba)$.
Working as in the proof of \cite[Proposition 3.2]{ASa}, $\GK\toba(W)<\infty$ where $W$ is of diagonal type with matrix $\bp=(p_{rs})_{1\le r,s\le 3}$. By direct computation, $p_{33}=1$ and $\widetilde{p}_{23}=\qti_{ik}\ne 1$, so we get a contradiction with \cite[Proposition 4.16]{AAH-memoirs}. Hence $x_{ij}=0$ in $\toba$.
\epf

Assume that $\I^+=\I^-=\emptyset$ and each $C_{\ell}$ is not of Cartan type $A_n$, $D_n$ with $q=-1$ (thus we know the eminent pre-Nichols algebra of each $\bq^{(\ell)}$) and $\GK\toba<\infty$. 
Then $\GK \toba^{(\ell)}<\infty$, so $\toba^{(\ell)}$ is a quotient of $\htoba_{\bq^{(\ell)}}$. Set
\begin{align}\label{eq:defn-htoba-non-connected}
\htoba_{\bq} &:= \bigotimes_{1\le \ell\le d} \htoba_{\bq^{(\ell)}}, & \hZ_{\bq} &:= \bigotimes_{1\le \ell\le d} \hZ_{\bq^{(\ell)}}, & 
\hOc{\bq} &=\bigcup_{1\le\ell\le d} \hOc{\bq^{(\ell)}},
\end{align}
extending the definitions we have made from the connected to the non-connected case. Hence we have an extension of $\N_{0}$-graded Hopf algebras
$\hZ_{\bq}\hookrightarrow \htoba_{\bq} \twoheadrightarrow \toba_{\bq}$, $\hZ_{\bq}$ is a polynomial ring in variables $z_{\ubeta}$, $\beta\in\hOc{\bq}$ of degree $\ubeta$,
and we may wonder if $\htoba_{\bq}$ is an eminent pre-Nichols algebra of $\bq$. We will see that this is the case, and then prove that the poset $\prefdgr{\bq}$ splits as the product of the posets $\prefdgr{\bq^{(\ell)}}$.

\begin{theorem}\label{thm:disconnected}
Let $\bq$ be such that $\I^+=\I^-=\emptyset$.
\begin{enumerate}[leftmargin=*,label=\rm{(\roman*)}]
\item\label{item:disconnected-pre-Nichols} If $\toba$ is an $\N_0$-graded pre-Nichols algebra of $\bq$ such that $\GK\toba<\infty$, then
\begin{align*}
\toba \simeq \otimes_{\ell=1}^d \toba^{(\ell)},
\end{align*}
and $\htoba_{\bq}$ is the eminent pre-Nichols algebra of $\bq$.

\item\label{item:disconnected-poset} There exists an anti-isomorphism of posets 
\begin{align*}
\prefdgr{\bq} \simeq \prod_{\ell=1}^{d} \prefdgr{\bq^{(\ell)}} \simeq \comp{\hOc{\bq}}.
\end{align*}
\end{enumerate}
\end{theorem}
\pf
\ref{item:disconnected-pre-Nichols} Notice that $\toba'=\otimes_{\ell=1}^d \toba^{(\ell)}$ is an $\N_0^{\theta}$-graded Hopf algebra with defining relations those defining each $\toba^{(\ell)}$ together with $x_{ij}=0$ for $i\in C_{k}$, $j\in C_{\ell}$, $k\ne \ell$. Thus Lemma \ref{lem:diff-conn-components} says that there exists a surjective map $\toba' \twoheadrightarrow \toba$ of $\N_0$-graded Hopf algebras, which is the identity on $V$, i.e. a map of pre-Nichols algebras of $\bq$. As we also have a map 
$\htoba_{\bq} \twoheadrightarrow \toba'$, the composition of both map gives a map of pre-Nichols algebras. As $\toba$ is arbitrary with finite $\GK$, $\htoba_{\bq}$ is eminent.

Hence we move to the description of $\N_0^{\theta}$-graded quotients of $\htoba_{\bq}$.
As before, we can use Proposition \ref{prop:twist-equiv-poset-preNichols} and Lemma \ref{lem:twist-equivalent-central} to assume that $\hZ_{\bq}$ is central. By Proposition \ref{prop:HopfidealB-HopfidealZ} we have to compute all the Hopf ideals of $\hZ_{\bq}$, which in turn (by taking graded duals) are classified by $\N_0^{\theta}$-graded Lie subalgebras of the nilpotent Lie algebra $\n_{\bq}=\prod_{1\le \ell\le d} \n_{\bq^{(\ell)}}$: any of these Lie subalgebras of $\n_{\bq}$ is the product of $\N_0^{\theta}$-graded Lie subalgebras of each $\n_{\bq^{(\ell)}}$, and these are classified by subsets $B^{(\ell)}$ of $\varDelta_{+}^{\bq^{(\ell)}}$ closed by sums, applying Theorem \ref{thm:Z-Delta-poset-ideals-sets} and Proposition \ref{prop:eminent-non-dist-poset} (depending on $\ell$). 

Coming back to $\hZ_{\bq}$, we quotient this Hopf algebra by $z_{\ubeta}$ for those $\ubeta\notin B:=\cup_{1\le \ell\le d} B^{(\ell)}$. We obtain then the pre-Nichols algebra
$$ \toba(\bq,B):=\htoba_{\bq}/\langle z_{\ubeta} | \ubeta\notin B \rangle \simeq \bigotimes_{1\le \ell\le d} \toba(q^{(\ell)},B^{(\ell)}), $$
and any $\toba\in\prefdgr{\bq}$ is of this shape, with $\toba^{(\ell)}\simeq \toba(q^{(\ell)},B^{(\ell)})$, $1\le \ell\le d$.

Now \ref{item:disconnected-poset} follows from \ref{item:disconnected-pre-Nichols} and the fact that a subset $B\in\hOc{\bq}$ is closed by sums if and only if each $B^{(\ell)}=B\cap\hOc{\bq^{(\ell)}}$ is closed by sum, since the sum $\alpha+\beta$ of two roots $\alpha\in\hOc{\bq^{(k)}}$, $\beta\in\hOc{\bq^{(\ell)}}$, $k\ne \ell$, is not a root.
\epf

\begin{remark}
Putting together Theorems \ref{thm:disconnected}, \ref{thm:connected-isom-posets} and \cite{C} we get a full description of the poset of all $\N_0^{\theta}$-graded pre-Nichols algebras with finite $\GK$ when all connected components are not points with label $\pm1$ neither of Cartan type $A_n$, $D_n$ and $q=-1$.
\end{remark}

\begin{remark}
Pre-Nichols algebras with finite $\GK$ for $\bq=\bq^{+}$ were classified in \cite[\S 3.4]{A-inf-dim-HA}. We do not include them in Theorem \ref{thm:disconnected} since we do not know at the moment how to control the interaction between $\toba^{+}$ and the pre-Nichols algebras $\toba^{-}$, $\toba^{(\ell)}$. 

At the same time, we do not have at the moment a description of all pre-Nichols algebras with finite $\GK$ for $\bq=\bq^{-}$ and for connected components of Cartan type $A_{\theta}$, $D_{\theta}$. Anyway, by Lemma \ref{lem:diff-conn-components} we may wonder that the poset of $\bq$ decomposes as the product of the posets of pre-Nichols of different connected components and that of $\bq^{-}$.
\end{remark}

\end{document}